\documentclass[a4paper,11pt]{article}

\usepackage{amsmath,amsthm,amstext,amscd,amssymb,euscript,mathrsfs}
\usepackage{calrsfs}
\usepackage{epsfig}

\usepackage{color}
\usepackage{esint}
\usepackage{pdfsync}
\usepackage{graphicx}
\usepackage{todonotes}

\newcommand{\epsi}{\epsilon}
\newcommand{\Z}{\mathbb Z}
\newcommand{\R}{\mathbb R}
\newcommand{\N}{\mathbb N}

\newcommand{\E}{\mathbb E}

\renewcommand{\phi}{\varphi}
\newcommand{\eps}{\ensuremath{\epsilon}}

\newcommand{\la}{\ensuremath{\Lambda}}
\newcommand{\si}{\ensuremath{\sigma}}

\numberwithin{equation}{section}

\def\1{{\mathchoice {\rm 1\mskip-4mu l} {\rm 1\mskip-4mu l}
{\rm 1\mskip-4.5mu l} {\rm 1\mskip-5mu l}}}

\newtheorem{theorem}{{\small T}{\scriptsize HEOREM}}[section]
\newtheorem{corollary}{{\bf{\small C}{\scriptsize OROLLARY}}}[section]
\newtheorem{proposition}{{\bf{\small P}{\scriptsize ROPOSITION}}}[section]
\newtheorem{lemma}{{\bf{\small L}{\scriptsize EMMA}}}[section]
\newtheorem{remark}{{\bf{\small R}{\scriptsize EMARK}}}[section]
\newtheorem{definition}{{\bf{\small D}{\scriptsize EFINITION}}}[section]

\renewenvironment{proof}[1]
{\noindent{{\bf{\small{P}{\scriptsize ROOF}}}.}\hspace{0.1cm} #1} {$\;\qed$\newline}

\newcommand{\beq}{\begin{eqnarray}}
\newcommand{\eeq}{\end{eqnarray}}

\newcommand{\ba}{\begin{align*}}
\newcommand{\ea}{\end{align*}}

\newcommand{\be}{\begin{equation}}
\newcommand{\ee}{\end{equation}}

\newcommand{\bl}{\begin{lemma}}
\newcommand{\el}{\end{lemma}}

\newcommand{\br}{\begin{remark}}
\newcommand{\er}{\end{remark}}

\newcommand{\bt}{\begin{theorem}}
\newcommand{\et}{\end{theorem}}

\newcommand{\bd}{\begin{definition}}
\newcommand{\ed}{\end{definition}}

\newcommand{\bp}{\begin{proposition}}
\newcommand{\ep}{\end{proposition}}

\newcommand{\bc}{\begin{corollary}}
\newcommand{\ec}{\end{corollary}}

\newcommand{\bpr}{\begin{proof}}
\newcommand{\epr}{\end{proof}}

\newcommand{\bi}{\begin{itemize}}
\newcommand{\ei}{\end{itemize}}

\newcommand{\ben}{\begin{enumerate}}
\newcommand{\een}{\end{enumerate}}

\newcommand{\caL}{{\mathcal L}}

\newcommand{\caU}{{\mathcal U}}

\newcommand{\tor}{\mathbb{T}}

\renewcommand{\(}{\left(}        \renewcommand{\)}{\right)}
\renewcommand{\[}{\left[}        \renewcommand{\]}{\right]}

     \newcommand{\nn}{\nonumber}

\newcommand{\red}[1]{\textcolor[rgb]{0,0,0}{#1}}

\begin{document}
\title{Asymmetric stochastic transport models \\
with ${\mathcal{U}}_q(\mathfrak{su}(1,1))$ symmetry}
\author{
Gioia Carinci$^{\textup{{\tiny(a)}}}$,
Cristian Giardin{\`a}$^{\textup{{\tiny(a)}}}$,\\
Frank Redig$^{\textup{{\tiny(b)}}}$,
Tomohiro Sasamoto $^{\textup{{\tiny(c)}}}$.\\\\
{\small $^{\textup{(a)}}$ Department of Mathematics, University of Modena and Reggio Emilia}\\
{\small via G. Campi 213/b, 41125 Modena, Italy}
\\
{\small $^{\textup{(b)}}$ Delft Institute of Applied Mathematics, Technische Universiteit Delft}\\
{\small Mekelweg 4, 2628 CD Delft, The Netherlands}\\
{\small $^{\textup{(c)}}$ Department of Physics, Tokyo Institute of Technology,}\\
{\small 2-12-1 Ookayama, Meguro-ku, Tokyo, 152-8550, Japan}\\
}
\maketitle
\begin{abstract}
By using the algebraic construction outlined in \cite{CGRS},
we introduce several Markov processes related to the
${\mathcal{U}}_q(\mathfrak{su}(1,1))$ quantum Lie algebra.
These processes serve as asymmetric transport models
and their
algebraic structure
easily allows to deduce
duality properties of the systems.
The results include:
(a) the asymmetric version of the Inclusion Process,
which is self-dual;
(b) the diffusion limit of this process, which is a natural
asymmetric analogue of the Brownian Energy Process
and which turns out to have the symmetric Inclusion
Process  as a dual process;
(c) the asymmetric analogue of the KMP
Process, which also turns
out to have a symmetric dual process.
We  give applications of the various duality relations
by computing exponential moments of the current.

%
%
\end{abstract}
\newpage
\tableofcontents
\newpage
\section{Introduction}

\subsection{Motivations}
Exactly solvable stochastic systems out-of-equilibrium have received considerable
attention in recent days 
\cite{Schuetz2000, GM, Derrida2007, QS2015p,  corwin2012kardar, BP}.
Often in the analysis of these models duality (or self-duality) is a crucial ingredient
by which the study of $n$-point correlations is reduced to the study of
$n$ dual particles. For instance, the exact current statistics
in the case of the asymmetric exclusion process is obtained by
solving the dual particle dynamics via Bethe ansatz \cite{Schutz, IS, borodin}.

The duality property has algebraic roots,
as was first noticed by Sch\"utz and Sandow for symmetric exclusion processes \cite{SS},
which is related to the  classical Lie algebra $\mathfrak{su}(2)$.
Next this symmetry approach was extended by  Sch\"utz  \cite{Schutz}
to the quantum Lie algebra ${\mathcal{U}}_q(\mathfrak{su}(2))$
in a representattion of spin $1/2$, thus providing self-duality
of the asymmetric exclusion process.
Recently Markov processes with the ${\mathcal{U}}_q(\mathfrak{su}(2))$
algebraic structure for higher spin value have been introduced and studied in \cite{CGRS}.
This lead to a family of non-integrable asymmetric generalization of the partial
exclusion process (see also \cite{M}).

In \cite{GKR, GKRV} the algebraic approach to duality
has been extended by connecting
duality functions to the algebra of operators commuting with
the  generator of the process.  In particular for the models of heat conduction
studied in  \cite{GKRV} the underlying algebraic structure
turned out to be ${\mathcal{U}}(\mathfrak{su}(1,1))$.
This class is richer than its fermionic counterpart
related to the classical Lie algebra ${\mathcal{U}}(\mathfrak{su}(2))$
which is at the root of processes of exclusion type.
In particular, the classical Lie algebra ${\mathcal{U}}(\mathfrak{su}(1,1))$
has been shown to be related to a large class of symmetric processes,
including: (a) an interacting particle system with attractive interactions
(inclusion process \cite{GKRV, GRV}); (b) interacting diffusion processes for heat
conduction (Brownian energy process \cite{GKRV, CGGR}); (c) redistribution models
of KMP-type \cite{KMP, cggr}. The dua\-li\-ties and self-dualities of all these processes arise naturally
from the symmetries which are built in the construction.

It is the aim of this paper to provide the  asymmetric version
of these models with (self)-duality property,
via the study of the deformed quantum Lie
algebra   ${\mathcal{U}}_q(\mathfrak{su}(1,1))$.
This provides a new class of bulk-driven non-equilibrium systems
with duality, which includes in particular an asymmetric version of
the KMP model  \cite{KMP}.
The diversity of models related to the classical
${\mathcal{U}}(\mathfrak{su}(1,1))$ will also appear here
in the asymmetric context where we consider the
quantum Lie algebra  ${\mathcal{U}}_q(\mathfrak{su}(1,1))$.

\subsection{Models and abbreviations}

For the sake of simplicity, we will use the following acronyms
in order to describe the class of new processes that arise
from our construction.

\begin{itemize}

\item[(a)] Discrete representations will provide  interacting particle systems in the class of {\em Inclusion Processes}.
For a parameter $k\in\mathbb{R}_+$, the Symmetric Inclusion Process version is denoted by
SIP$(k)$, and  {\bf ASIP$(q,k)$} is the corresponding asymmetric version, with asymmetry parameter $q \in (0,1)$.
\item[(b)] Continuous representations give rise to diffusion processes in the class of {\em Brownian  Energy Processes}. For  $k\in\mathbb{R}_+$, the Symmetric Brownian Energy Process is denoted by  BEP$(k)$, and {\bf ABEP$(\sigma,k)$} is the asymmetric version with asymmetry parameter $\sigma>0$.
\item[(c)] By instantaneous thermalization, redistribution models are obtained,  where energy or particles are redistributed at Poisson event times. This class includes the thermalized version of ABEP$(\sigma,k)$, which is denoted by Th-ABEP$(\sigma,k)$. In the particular case $k=1/2$ the Th-ABEP$(\sigma,k)$ is called the asymmetric KMP (Kipnis-Marchioro-Presutti) model, denoted by {\bf AKMP$(\sigma)$}, which becomes the 
KMP model as $\sigma \to 0$. 
The instantaneous thermalization of the ASIP$(q,k)$ yields the Th-ASIP$(q,k)$ process.
\end{itemize}

\subsection{Markov processes with algebraic structure}

In \cite{CGRS} we constructed a generalization of the asymmetric exclusion process,
allowing $2j$ particles per site with self-duality properties reminiscent
of the self-duality of the standard ASEP found initially by Sch\"{u}tz in \cite{Schutz}.
This construction followed a general scheme where one starts from
the Casimir operator $C$ of the quantum Lie algebra ${\mathcal{U}}_q(\mathfrak{su}(2))$, and applies a coproduct to
obtain an Hamiltonian $H_{i,i+1}$ working on the occupation number
variables 
at sites $i$ and $i+1$.
The operator $H=\sum_{i=1}^L H_{i,i+1}$ then naturally allows a rich class of
commuting operators (symmetries), obtained from the $n$-fold coproduct applied
to any generator of the algebra. This  operator $H$ is
not yet the generator of a Markov process.
But $H$ allows a strictly  positive ground state, which can also be constructed
from the symmetries applied to a trivial ground state.
Via a ground state transformation, $H$ can then be turned into a Markov generator $L$
of a jump process where particles hop between nearest neighbor sites and at most
$2j$ particles per site are allowed. The symmetries of $H$ directly translate into the
symmetries of $L$, which in turn directly translate into self-duality functions.

This construction is in principle applicable to every quantum Lie algebra
with a non-trivial center.
However, it is not guaranteed that a Markov generator can be obtained. This
depends on the chosen representation of the generators of the algebra, and the choice
of the co-product.
Recently the construction has been applied to algebras with higher rank,
such as ${\mathcal{U}}_q(\mathfrak{gl}(3))$ \cite{bel,kuan}
or ${\mathcal{U}}_q(\mathfrak{sp}(4))$ \cite{kuan},
yielding  two-component  asymmetric exclusion process with
multiple conserved species of particles.


\subsection{Informal description of main results}

In \cite{GKRV} we introduced a class of processes with $\mathfrak{su}(1,1)$ symmetry
which in fact arise from this construction for the Lie algebra $\caU (\mathfrak{su}(1,1))$.
In this paper we look for natural asymmetric versions of the processes
constructed in \cite{GKRV}, and \cite{cggr}. In particular the natural
asymmetric analogue of the KMP process is a target.
The main results are the following
\begin{itemize}
\item[(a)] {\it Self-duality of ASIP$(q,k)$.}
We proceed
via the same construction as in \cite{CGRS} for the algebra ${\mathcal{U}}_q(\mathfrak{su}(1,1))$ to find the
ASIP$(q,k)$ which is the ``correct'' asymmetric analogue of
the SIP$(k)$. The parameter $q$ tunes the asymmetry: $q\to 1$ gives back
the SIP$(k)$.
This process is then via its construction self-dual with a non-local
self-duality function.
\item[(b)]  {\it Duality between ABEP$(\si,k)$ and SIP$(k)$.}
We then show that in the limit $\epsi\to 0$ where simultaneously the asymmetry is
going to zero ($q=1-\epsi \si$ tends to unity), and the number of particles to infinity 
$\eta_i= \lfloor\epsi^{-1} x_i\rfloor$, we obtain
a diffusion process ABEP$(\si,k)$ which is reminiscent of the Wright-Fisher diffusion with
mutation and a selective drift. As a consequence of self-duality of ASIP$(q,k)$ we
show that this diffusion process is dual to the SIP$(k)$, i.e., the dual process is
symmetric, and the asymmetry is in the duality function.
Notice that this is the first example of duality between a truly asymmetric
system (i.e. bulk-driven)  and a symmetric system (with zero current).

\item[(c)]  {\it Duality of instantaneous thermalization models.}
Finally, we then consider instantaneous thermalization of ABEP$(\si,k)$ to
obtain an asymmetric energy redistribution model of KMP type. Its dual
is the instantaneous thermalization of the SIP$(k)$ which for $k=1/2$ is exactly the dual KMP
process.
\end{itemize}

\subsection{Organization of the paper}

The rest of our paper is organized as follows. In section \ref{asipabep} we introduce
the process ASIP$(q,k)$. After discussing some limiting cases,
we show that this process has reversible profile product measures
on $\mathbb{Z}_+$ (but not on $\mathbb{Z}$).

In section \ref{abep-sec} we consider the  weak asymmetry limit of ASIP$(q,k)$.
This leads to the diffusion process ABEP$(\si,k)$, that also has
reversible inhomogeneous product measures on the half-line.
We prove that ABEP$(\sigma,k)$ is a genuine non-equilibrium asymmetric system
in the sense that it has a non-zero average current.
Nevertheless in the last part of section \ref{abep-sec} we show that
the ABEP$(\sigma,k)$ can be mapped -- via a global change of coordinates --
to the BEP$(k)$, which is a symmetric system with zero-current.
In section  \ref{otto} this is also explained in the framework of the representation theory of the classical Lie algebra
${\mathcal{U}}(\mathfrak{su}(1,1))$.

In section  \ref{akmp-section} we introduce the instantaneous thermalization limits
of both ASIP$(q,k)$ and ABEP$(\si,j)$ which are a particle, resp.\ energy,
redistribution model at Poisson event times. This provides asymmetric
redistribution models of KMP type.

In section \ref{dualsect} we introduce the self-duality of the ASIP$(q,k)$ and prove
various other duality relations that follow from it. In particular,
once the self-duality of ASIP$(q,k)$ is obtained, duality of
ABEP$(\sigma,k)$ with SIP$(k)$ follows from a limiting procedure
which is proved in Section \ref{subsect52}.
In the limit of an infinite number of particles with weak-asymmetry,
the original process scales to ABEP$(\sigma,k)$,
whereas in the dual process the asymmetry disappears because
 the number of particles is finite.
 Next the self-duality and duality of thermalized models
is derived in Section \ref{subsect53}.

In section \ref{computations} we illustrate the use of the
duality relations in various computations of exponential  moments of currents.
Finally, the last section is devoted to the full construction of the ASIP$(q,k)$ from
a ${\mathcal{U}}_q(\mathfrak{su}(1,1))$ symmetric quantum Hamiltonian and the
proof of self-duality from the symmetries of this Hamiltonian.

\section{The Asymmetric Inclusion Process ASIP$(q,k)$}\label{asipabep}
\subsection{Basic notation}
We will consider as underlying lattice the finite lattice $\la_L=\{1,\ldots,L\}$
or the periodic lattice $\tor_L=\Z/L\Z$.
At the sites of $\la_L$ we allow an arbitrary number of particles.
The particle system configuration space is $\Omega_L=\N^{\la_L}$. Elements of $\Omega_L$
are denoted by $\eta, \xi$ and for $\eta\in\Omega_L$, $i \in \Lambda_L$, we denote by $\eta_i\in \N$ the number
of particles at site $i$.
For $\eta\in \Omega_L$ and $i,j\in \la_L$ such that $\eta_i>0$,
we denote by $\eta^{i,j}$ the configuration
obtained from $\eta$ by removing one particle from $i$ and
putting it at $j$.\\

We need some further notation of $q$-numbers.
For $q\in(0,1)$  and $n\in \mathbb{N}_0$ we introduce the {$q$-number}
\begin{equation}
\label{q-num}
[n]_q=\frac{q^n-q^{-n}}{q-q^{-1}}
\end{equation}
satisfying the property $\lim_{q\to1} [n]_q = n$.
The first $q$-number's are thus given by
$$[
0]_q = 0, \quad\quad\quad [1]_q=1, \quad \quad \quad [2]_q=q+ q^{-1}, \quad \quad \quad [3]_q=q^2+1+q^{-2}, \quad \dots
$$
We also introduce the $q$-factorial
$$
[n]_q!:=[n]_q \cdot [n-1]_q \cdot \dots \cdot [1]_q \;,
$$
and the $q$-binomial coefficient
$$
\binom{n}{m}_q:=\frac{[n]_q!}{[m]_q! [n-m]_q!}  \;.
$$
Further we denote
\be\label{Poch}
(a;q)_{m} := (1-a)(1-aq)\cdots (1-aq^{m-1})\;.
\ee

\subsection{The ASIP$(q,k)$ process}
We introduce the process in finite volume by specifying its generator.

\bd[ASIP(q,k) process]
\label{def}
\vskip.2cm
\noindent
\begin{enumerate}
\item The ASIP$(q,k)$ with closed boundary conditions is defined
as the Markov process on $\Omega_L$ with generator
defined on functions $f:\Omega_L\to\R$
\begin{eqnarray}
\label{gen}
&&({\cal L}^{ASIP(q,k)}_{(L)}f)(\eta):=\sum_{i=1}^{{L-1}} ({\cal L}^{ASIP(q,k)}_{i,i+1}f)(\eta) \qquad \text{with} \nonumber\\
({\cal L}^{ASIP(q,k)}_{i,i+1}f)(\eta)
& := & q^{\eta_i-\eta_{i+1}+(2k-1)} [\eta_i]_q [2k+\eta_{i+1}]_q (f(\eta^{i,i+1}) - f(\eta)) \nonumber \\
& +  & q^{\eta_i-\eta_{i+1}-(2k-1)} [2k+\eta_i]_q [\eta_{i+1}]_q (f(\eta^{i+1,i}) - f(\eta))
\end{eqnarray}
\item The ASIP$(q,k)$ with periodic boundary conditions is defined
as the Markov process on $\N^{\tor_L}$ with generator
\be
({\cal L}^{ASIP(q,k)}_{(\tor_L)}f)(\eta) :=\sum_{i\in \tor_L} ({\cal L}^{ASIP(q,k)}_{i,i+1} f)(\eta)
\ee
\end{enumerate}
\ed
\noindent
Since in finite volume we always start with finitely many particles, and the total particle number is conserved, the process is automatically well defined as a finite state space continuous time Markov chain. Later on (see Section \ref{Infi})  we will consider expectations of the self-duality functions in the  infinite volume limit. In this way we can deal with relevant infinite volume expectations without having to solve the full existence problem of the ASIP$(q,k)$ in infinite volume for a generic initial data. This might actually be an hard problem due to the lack of monotonicity.
\subsection{Limiting cases}
The ASIP $(q,k)$ degenerates to well known interacting particle systems when its parameters take the limiting values $q \to 1$ and $k \to \infty$ recovering the cases of symmetric evolution or totally asymmetric zero range interaction. Notice in particular that these two limits do not commute.
\begin{itemize}
\item{Convergence to symmetric processes}
\begin{itemize}
\item[i) ] {\bf $\mathbf{q \to 1, k}$ fixed:}
The ASIP$(q,k)$ reduces to the SIP$(k)$, i.e. the
 Symmetric Inclusion Process with parameter $k$.
All the results of the present paper  apply also to this symmetric case.
In particular, in the limit $q\to 1$, the self-duality functions that will be given in theorem
\ref{mainself} below converge to the self-duality functions of the SIP$(k)$ (given in \cite{cggr}).

\item[ii) ]{ $\mathbf{q \to 1, k \to \infty}${\bf :}} Furthermore, when the symmetric  inclusion process is time changed so that time is scaled down
by a factor $1/2k$, then in the limit $k\to\infty$ the symmetric inclusion converges weakly
in path space to a system of symmetric independent random walkers (moving at rate 1).
\end{itemize}
\item Convergence to totally asymmetric processes
\begin{itemize}
\item[ iii)]{\bf $\mathbf{k \to \infty, q}$ fixed:}
If the limit $k\to\infty$ is performed first, then a totally asymmetric system is obtained under proper time rescaling. Indeed, by multiplying the ASIP$(q,k)$ generator by $(1-q^2)q^{4k-1}$ one has
\begin{eqnarray*}
(1-q^2)q^{4k-1} \[\mathcal L^{ASIP}_{i,i+1}f\](\eta)&=&q^{4k}\; \frac{(q^{2\eta_i}-1)(q^{4k}-q^{-2\eta_{i+1}})}{(1-q^2)} \[f(\eta^{i,i+1})-f(\eta))\] \nn \\
&+&\; \frac{(q^{-2\eta_{i+1}}-1)(1-q^{2\eta_i+4k})}{(q^{-2}-1)} \[f(\eta^{i+1,i})-f(\eta))\] \nn
\end{eqnarray*}
Therefore, considering  the family of processes $y^{(k)}(t):=\{y_i^{(k)}(t)\}_{i \in \Lambda_L}$ labeled by $k\ge0$ and
defining
$$
y_i^{(k)}(t) :=\eta_i((1-q^2)q^{4k-1}t)
$$
one finds that in the limit
$k\to \infty$ the process $y^{(k)}(t)$ converges weakly to  the Totally Asymmetric Zero Range process $y(t)$ with
generator  given by:
\be\label{genT11}
({\cal L}^{q-\text{TAZRP}}_{su(1,1)}f)(y) =\sum_{i=1}^{L-1}\,\frac{q^{-2y_{i+1}} -1}{q^{-2}-1} \;[f(y^{i+1,i})-f(y)], \qquad f: \Omega_L\to \mathbb R
\ee
In this system,  particles jump to the left only with rates that are monotone increasing functions of the occupation variable of the departure site.
Note that the rates are unbounded for $y_{i+1}\to\infty$, nevertheless the process is well defined even in the infinite volume,
as it belongs to the class considered in \cite{seppalainen}.
This is to be compared to the case of the deformed algebra $U_q(\mathfrak{sl}_2)$ \cite{CGRS} whose scaling limit with infinite spin
is given by \cite{borodin}
\be\label{genT2}
({\cal L}^{(q-\text{TAZRP})}_{su(2)}f)(y) =\sum_{i=1}^{L-1} \,\frac{1-q^{2y_i}}{1-q^2} \;[f(y^{i,i+1})-f(y)],
 \qquad f: \Omega_L \to \mathbb R
\ee
Here particles jump to the right only  with rates that are also a monotonous increasing function
of the occupation variable of the departure site, however now it is a bounded function
approaching $1$ in the limit $y_{i}\to\infty$.
In \cite{FV} it is proved that the totally asymmetric zero range process \eqref{genT2}
is in the KPZ universality class. It is an interesting open problem to prove
or disprove that the same conclusion holds true for  \eqref{genT11} \cite{povo}.
We remark that the rates of \eqref{genT11} are (discrete) convex function and
this also translates into convexity of the stationary
current $j(\rho)$ as a function of the density $\rho$,
whereas for \eqref{genT2} we have concave relations.
\item[iv)] {\bf $\mathbf{k \to \infty, q \to 1}$:} In the limit $q\to 1$ the
zero range process in \eqref{genT11} reduces to a system of totally
asymmetric independent walkers. This is to be compared to item ii)
where symmetric walkers were found if the two limits
were performed in the reversed order.
\end{itemize}
\end{itemize}
\subsection{Reversible profile product measures}
Here we describe the reversible measures of ASIP$(q,k)$.
\begin{theorem}[{Reversible measures of ASIP$(q,k)$}]
\label{basicproptheorem1}
For all $L\in\mathbb{N}, L\geq 2$, the following results hold true:
\begin{enumerate}
\item[1.)]
the ASIP$(q,k)$ on $\Lambda_L$ with closed boundary conditions
admits a family labeled by $\alpha$ of reversible product measures with marginals given by
\be
\label{stat-meas}
\mathbb{P^{(\alpha)}}(\eta_i = n) = \frac{\alpha^n}{Z^{(\alpha)}_i} \,{\binom{n+2k-1}{n}_q} \cdot  q^{4kin}  \qquad\qquad n \in \mathbb N
\ee
 for $i \in \Lambda_L$ and $\alpha \in [0, q^{-(2k+1)})$ (with the convention $\binom{2k-1}{0}_q = 1$).
 The normalization is
\be \label{Z}
Z_i^{{(\alpha)}} = \sum_{n=0}^{+\infty}  {\binom{n+2k-1}{n}_q} \cdot  \alpha^n q^{4kin} = \frac{1}{(\alpha q^{4ki-(2k-1)};q^2)_{2k}}
\ee
and for this measure
\be
\label{somma}
\mathbb{E}^{(\alpha)}(\eta_i) = \sum_{l=0}^{2k-1} \frac{1}{q^{-2l}(\alpha q^{4ki - 2k +1})^{-1} -1}\;.
\ee
\item[2.)]
The ASIP$(q,k)$ process on the torus $\mathbb{T}_L$ with periodic boundary condition does not admit
homogeneous product measures.
\end{enumerate}
\end{theorem}

\bpr
The proof of item 2.) is similar to the proof of Theorem 3.1, item d) in \cite{CGRS} and
we refer the reader to that paper for all details.
To prove item 1.) consider the detailed balance relation
\be\label{DB}
\mu(\eta)c_q(\eta, \eta^{i,i+1})=\mu(\eta^{i,i+1})c_q(\eta^{i,i+1},\eta)
\ee
where the hopping rates are given by
\begin{equation*}
c_q(\eta, \eta^{i,i+1})= q^{\eta_i-\eta_{i+1}+2k-1} [\eta_i]_q [2k+\eta_{i+1}]_q
\end{equation*}
\begin{equation*}
c_q(\eta^{i,i+1},\eta)= q^{\eta_i-\eta_{i+1}-2k-1} [2k+\eta_i-1]_q [\eta_{i+1}+1]_q
\end{equation*}
and  $\mu$ denotes a reversible measure.
Suppose now that $\mu$ is a product measure of the form $\mu=\otimes_{i=1}^L\mu_i$.
Then  \eqref{DB} holds if and only if
\be
\mu_i(\eta_i-1)\mu_{i+1}(\eta_{i+1}+1) q^{-2k} [2k+\eta_i-1]_q [\eta_{i+1}+1]_q = \mu_i(\eta_i)\mu_{i+1}(\eta_{i+1}) q^{2k} [\eta_i]_q [2k+\eta_{i+1}]_q
\ee
which implies that there exists $\alpha \in \mathbb R$ so that for all $i\in \Lambda_L$
\be \label{muuu}
\frac{\mu_i(n)}{\mu_i(n-1)}= \alpha q^{4ki} \frac{[2k+n-1]_q}{[n]_q}\;.
\ee
Then \eqref{stat-meas}  follows from \eqref{muuu} after using an induction argument on $n$.
The normalization $Z_i^{{(\alpha)}}$ is computed by using Corollary 10.2.2 of \cite{AAR}.
We have that
\be
Z_i^{{(\alpha)}}<\infty \quad \text{if and only if } \quad 0 \le \alpha < q^{-4ki+(2k-1)} \quad \text{for any}\: i \in \Lambda_L
\ee
As a consequence (since $q<1$ and $i=1$ is the worst case) $\alpha$ must belong to the interval $[0,q^{-(2k+1)})$.
The expectation \eqref{somma} is obtained by exploiting the identity
$$
\mathbb{E}^{(\alpha)}(\eta_i) =  \alpha \frac{d}{d\alpha} \log Z_i^{{(\alpha)}}.
$$
\epr
The following comments are in order:
\begin{itemize}
\item[i)]
{\em vanishing asymmetry}: in the limit $q \to 1$ the reversible product measure of ASIP$(q,k)$
converges to a product of Negative Binomial distributions with shape parameter $2k$ and
success probability $\alpha$, which are the reversible measures of
the SIP$(k)$ \cite{cggr}.
\item[ii)]
 {\em monotonicity of the profile:}
the average occupation number $\mathbb{E}^{(\alpha)}(\eta_i)$ in formula \eqref{somma}
is a decreasing function of $i$, and $\lim_{i\to\infty} \mathbb{E}^{(\alpha)}(\eta_i) = 0$.
\item[iii)]
{\em infinite volume}:
the reversible product measures with marginal \eqref{stat-meas} are also well-defined in the limit
$L\to\infty$. One could go further to $[-M,\infty)\cap \mathbb{Z}$ for $\alpha < q^{4kM + 2k - 1}$
(but not to the full  line $\mathbb{Z}$). These infinite volume measure concentrate on
configurations with a finite number of particles, and thus are the analogue of the profile measures
in the asymmetric exclusion process \cite{Liggett}.
\end{itemize}

\section{The Asymmetric Brownian Energy Process ABEP$(\sigma,k)$}
\label{abep-sec}
Here we will take the limit of weak asymmetry $q= 1-\epsi\si\to 1$ ($\epsi\to 0$) combined with the number
of particles proportional to $\epsi^{-1}$, going to infinity, and work with rescaled particle numbers
$x_i= \lfloor \epsi\eta_i\rfloor$. Reminiscent of scaling limits in population dynamics, this leads to a diffusion process
of Wright-Fisher type \cite{CGGR}, with $\si$-dependent drift term, playing the role of a selective drift in the
population dynamics language, or bulk driving term in the non-equilibrium statistical physics language.

\subsection{Definition}
We  define the ABEP$(q,k)$ process via its generator.
It has state space ${\mathcal X}_L=(\mathbb R_+)^L$, $\mathbb R_+:=[0, +\infty)$. Configurations are denoted by   $x\in {\mathcal X}_L$, with
$x_i$ being interpreted as the energy at site $i\in \Lambda_L$.

\bd[ABEP$(\sigma,k)$ process]\label{def3}
\vskip.2cm
\noindent
\begin{enumerate}

\item
Let $\sigma>0$ and  $k\ge 0$.
 The Markov process ABEP$(\sigma,k)$ on the state space ${\mathcal X}_L$
 with
{closed} boundary conditions
is defined by the generator working on the core of smooth functions $f: {\mathcal X}_L\to \R$ via
\be
[{\cal L}^{ABEP^{(\sigma,k)}}_{(L)}f](x) =\sum_{i=1}^{{L-1}} [{\cal L}^{ABEP^{(\sigma,k)}}_{i,i+1}f](x)
\ee
with
\begin{eqnarray}
\label{genABEP}
&&  \[{\cal L}^{ABEP^{(\sigma,k)}}_{i,i+1}f\](x)
 = \frac 1{4\sigma^2}\, (1-e^{-2\sigma x_i})(e^{2\sigma x_{i+1}}-1)\(\frac{\partial}{\partial{x_i}}-\frac{\partial}{\partial{ x_{i+1}}}\)^2f(x)\nonumber \\
&&  -\frac{1}{2\sigma}\, \bigg\{ (1-e^{-2\sigma x_i})(e^{2\sigma x_{i+1}}-1)+2k\(2-e^{-2\sigma x_i} - e^{2\sigma x_{i+1}}\)\bigg\}\(\frac{\partial}{\partial{x_{i}}}-\frac{\partial}{\partial{x_{i+1}}}\)f(x)\nonumber
\end{eqnarray}
\item The ABEP$(\sigma,k)$ with periodic boundary conditions is defined
as the Markov process on $\mathbb R_+^{\tor_L}$ with generator
\be
[{\cal L}^{ABEP(\sigma,k)}_{(\tor_L)} f](x):=\sum_{i\in \tor_L} [{\cal L}^{ABEP(\sigma,k)}_{i,i+1}f](x)
\ee
\end{enumerate}
\ed

The ABEP$(\sigma,k)$ is a genuine asymmetric non-equilibrium system, in the sense that
its translation-invariant stationary state may sustain a non-zero current.
To see this, let $\mathbb{E}$ denote expectation with respect to the
translation invariant measure for  the ABEP$(\sigma,k)$ on  $\mathbb T_L$.
Let $f_i(x):=x_i$, then from \eqref{genABEP} we have
\begin{eqnarray}
 [{\cal L}^{ABEP^{(\sigma,k)}}f_i](x)= \Theta_{i,i+1}(x)- \Theta_{i-1,i}(x)
 \end{eqnarray}
 with
 \begin{eqnarray}
 \label{corri}
 \Theta_{i,i+1}(x)=   - \frac{1}{2\sigma}\, \bigg\{ (1-e^{-2\sigma x_i})(e^{2\sigma x_{i+1}}-1)+2k\(2-e^{-2\sigma x_i} - e^{2\sigma x_{i+1}}\)\bigg\}
\end{eqnarray}
So we have 
\begin{equation*}
\frac{d}{dt} \E_x \[f_i(x(t))\]= \E_x\[\Theta_{i, i+1} (x(t))\]-\E_x\[\Theta_{i-1,i}(x(t))\]
\end{equation*}
and then, {from the continuity equation we have that},  in a translation invariant state, $\mathcal{J}_{i,i+1}:=-\E\[\Theta_{i,i+1}\]$ is the {instantaneous} stationary current over the
edge $(i,i+1)$.
Thus we have  the following
\bp[Non-zero current of ABEP$(\sigma,k)$]
\label{non-zero}
$$
\mathcal{J}_{i,i+1} = - \mathbb E\[\Theta_{i,i+1}\] < 0  \qquad \text{if} \quad k>1/2
$$
and
$$
\mathcal{J}_{i,i+1} = - \mathbb E\[\Theta_{i,i+1}\] > 0  \qquad \text{if} \quad k=0\;.
$$
\ep
\bpr
In the case $k>1/2$, taking expectation of \eqref{corri} we obtain
$$
\mathbb E\[\Theta_{i,i+1}\] =
\frac{1}{2\sigma}\, \bigg\{ (1- 4k) + (2k - 1) \mathbb{E} (e^{2\sigma x_{i+1}} + e^{-2\sigma x_i})  + \mathbb{E}(e^{2\sigma(x_{i+1}-x_i)})\bigg\}
$$
Since expectation in the translation invariant stationary state of local variables are the same on each site
and $\cosh(x) \ge 1$ one obtains
$$
\mathbb E\[\Theta_{i,+1}\] \ge
\frac{1}{2\sigma}\,  \left\{(1-4k) + {2}(2k-1)  + \mathbb{E}\[e^{2\sigma(x_{i+1}-x_i)}\]\right\}
$$
Furthermore, Jensen inequality  and translation invariance  implies that
$$
\mathbb{E} \[\Theta_{i,i+1}\] >
\frac{1}{2\sigma}\, \Big\{ (1-4k) +{2} (2k-1)  + 1 \Big\} =0
$$
In the case $k=0$ one has
$$
\mathbb E\[\Theta_{i,i+1}\] =
\frac{1}{2\sigma}\, \mathbb E \Big[ (1-e^{-2\sigma x_i})(1- e^{2\sigma x_{i+1}})\Big] < 0
$$
 { which is negative because the function is negative a.s}.
\epr
\vskip.5cm
\subsection{Limiting cases}
\begin{itemize}
\item Symmetric processes
\begin{itemize}
\item[i)] {\bf $\mathbf{\sigma \to 0, k} $ fixed:} we recover  the Brownian Energy Process with
parameter k, BEP$(k)$ (see \cite{cggr}) whose generator is
\begin{eqnarray}
\label{genBEP}
{\cal L}^{BEP^{(k)}}_{i,i+1}
 = x_ix_{i+1}\(\frac{\partial}{\partial{x_i}}-\frac{\partial}{\partial{ x_{i+1}}}\)^2 -2k(x_i-x_{i+1})\(\frac{\partial}{\partial{x_{i}}}-\frac{\partial}{\partial{x_{i+1}}}\)
\end{eqnarray}
\item[ii)] {\bf $\mathbf{\sigma \to 0, k \to \infty} $:}
 under the time rescaling $t\to t/2k$, one finds
that in the limit $k\to\infty$ the BEP$(k)$ process scales to a
symmetric deterministic system evolving
with generator
\begin{eqnarray}\label{DEP}
\[{\cal L}^{DEP}_{i,i+1}f\](x)
& = &
- (x_i -  x_{i+1})  \(\frac{\partial}{\partial{x_{i}}}-\frac{\partial}{\partial{x_{i+1}}}\)f(x)
\end{eqnarray}
This deterministic system is symmetric in the sense that if the initial condition is given by
$(x_i(0),x_{i+1}(0))= (a,b)$ then the asymptotic  solution is
given by the fixed point $\left(\frac{a+b}{2},\frac{a+b}{2}\right)$
where the  initial total energy $a+b$  is equally shared among the two sites.

\end{itemize}
\item Wright-Fisher diffusion
\begin{itemize}
\item[iii)] {\bf $\mathbf{\sigma \simeq 0, k} $ fixed:}
the ABEP$(\sigma,k)$ on the simplex can be read as a Wright Fisher model with mutation and selection,
however we have not been able to find in the literature the specific form of selection
appearing in \eqref{genABEP} (see \cite{CGGR} for the analogous result when $\sigma=0$). To first order in $\sigma$ one recovers the standard
Wright-Fisher model with constant mutation $k$ and selection $\sigma$
\begin{eqnarray}
{\cal L}^{WF(\sigma,k)}
=  x_i x_{i+1} \(\frac{\partial}{\partial{x_i}}-\frac{\partial}{\partial{ x_{i+1}}}\)^2
 -
\left(2 \sigma x_i x_{i+1} + 2k(x_i -  x_{i+1}) \right)  \(\frac{\partial}{\partial{x_{i}}}-\frac{\partial}{\partial{x_{i+1}}}\)\nonumber
\end{eqnarray}
\end{itemize}

\item Asymmetric Deterministic System
\begin{itemize}
\item[iv)] {\bf $\mathbf{k \to \infty, \sigma} $ fixed:}
if the limit $k\to\infty$ is taken directly on the ABEP$(\sigma,k)$ then, by  time rescaling
$t\to t/2k$ one arrives at an asymmetric deterministic system with generator
\begin{eqnarray}
\label{genADEP}
 {\cal L}^{ADEP^{(\sigma)}}_{i,i+1}
& = &
- \frac{1}{2\sigma}\, \(2-e^{-2\sigma x_i} - e^{2\sigma x_{i+1}}\)\(\frac{\partial}{\partial{x_{i}}}-\frac{\partial}{\partial{x_{i+1}}}\)
\end{eqnarray}
This deterministic system is asymmetric in the sense that if the initial condition is given by
$(x_i(0),x_{i+1}(0))= (a,b)$ then the asymptotic  solution is
given by the fixed point
$$(A,B):=\left(\frac{1}{2\sigma}\ln\(\frac{1+e^{2\sigma(a+b)}}{2}\), a+b - \frac{1}{2\sigma}\ln\(\frac{1+e^{2\sigma(a+b)}}{2}\)\right)$$
where $A>B$.
\item[v)] {\bf $\mathbf{k \to \infty, \sigma \to 0} $:}  in the limit $\sigma\to 0$ \eqref{genADEP} converges to \eqref{DEP} and one
 recovers again the symmetric equi-distribution between the two sites of DEP process with generator \eqref{DEP}.

\item[vi)]{\bf $\mathbf{k \to \infty, \sigma \to \infty} $:}  in the limit $\sigma\to\infty$ one has the totally asymmetric stationary solution $(a+b,0)$.


\end{itemize}
\end{itemize}

\subsection{The ABEP$(\sigma,k)$ as a diffusion limit of  ASIP$(q,k)$.}

Here we show that the ABEP$(\sigma,k)$ arises from the ASIP$(q,k)$ in a limit of vanishing asymmetry and infinite
particle number.

\bt[Weak asymmetry limit of ASIP$(q,k)$]
\label{teo1}
Fix $T>0$.
Let $\{\eta^\epsi(t):0\leq T\}$ denote the
ASIP$(1-\si\epsi, k)$ starting from initial condition
$\eta^\epsi(0)$. Assume that
\be\label{initialco}
\lim_{\epsi\to  0} \epsi\eta^\epsi(0)= x\in {\mathcal X}_L
\ee
Then as $\epsi\to 0$, the process $\{\eta^\epsi(t):0\leq t\leq T\}$ converges weakly on
path space to the ABEP$(\si, k)$ starting from $x$.
\et
\bpr
The proof follows the lines of the corresponding results in population dynamics literature,
i.e., Taylor expansion of the generator and keeping the relevant orders.
Indeed, by the Trotter-Kurtz theorem \cite{Liggett}, we have to prove that on the core
of the generator of the limiting process, we have convergence of generators.
Because the generator is a sum of terms working on two variables, our theorem
follows from the  computational lemma below. \epr
\bl
If $\eta^\epsi\in \Omega_L$ is such that $\epsi\eta^\epsi\to x\in {\mathcal X}_L$ then, for
every smooth function $F: {\mathcal X}_L\to\R$, and for every $i\in \{1,\ldots, L-1\}$
we have
\be\label{geniconv}
\lim_{\epsi\to 0} ({\cal L}_{i,i+1}^{ASIP(1-\epsi \si, k)} F_\epsi)(\eta^\epsi)=
{\cal L}_{i,i+1}^{ABEP(\si, k)} F(x)
\ee
where $F_\epsi(\eta)= F(\epsi \eta)$, $\eta\in \Omega_L$.
\el
\bpr
Define $x^\epsi= \epsi \eta^\epsi$.
Then we have, by the regularity assumptions on $F$ that
\beq
&&F_\epsi ((\eta^\epsi)^{i,i+1})- F_\epsi(\eta)
\nonumber\\
&=&\eps\left (\frac{\partial }{\partial x_{i+1}}-\frac{\partial }{\partial x_{i}}\right )F(x^\eps)+\eps^2 \left (\frac{\partial }{\partial x_{i}}-\frac{\partial }{\partial x_{i+1}}\right )^2 F(x^\eps)+{  O(\eps^3)}
\eeq
and similarly
\beq
&&F_\epsi ((\eta^\epsi)^{i+1,i})- F_\epsi(\eta)\nonumber\\
&=&
 -\eps\left (\frac{\partial }{\partial x_{i+1}}-\frac{\partial }{\partial x_{i}}\right )F(x^\eps)+\eps^2 \left (\frac{\partial }{\partial x_{i}}-\frac{\partial }{\partial x_{i+1}}\right )^2 F(x^\eps)+{  O(\eps^3)}
\eeq
Then using
$q= 1-\epsi \si$, and
\be
(1-\epsi\si)^{x^\epsi_i/\epsi}= e^{-\si x_i} -2 x_i \si^2 e^{-2\si x_i}\epsi + O(\epsi^2)
\nonumber
\ee
straightforward computations give
\begin{eqnarray*}
\[{\cal L}_{i,i+1}^{\eps}F\](x^\eps)&=&\[B_\eps(x^\eps)\left(\frac{\partial}{\partial x_{i+1}}-\frac{\partial}{\partial x_{i}}\right) +
D_\eps(x^\eps)  \left ( \frac{\partial }{\partial x_{i}}-\frac{\partial }{\partial x_{i+1}}\right  )^2\]F(x^\eps) + O(\eps)
\end{eqnarray*}
with
\begin{eqnarray}
&& B_\eps(x)= \frac{1}{2\sigma}\, \bigg\{ (1-e^{-2\sigma x_i})(e^{2\sigma x_{i+1}}-1)+2k\(2-e^{-2\sigma x_i} - e^{2\sigma x_{i+1}}\)\bigg\} +  O(\eps)\nonumber \\
&& D_\eps(x)= \frac 1{4\sigma^2}\, (1-e^{-2\sigma x_i})(e^{2\sigma x_{i+1}}-1)+  O(\eps)
\end{eqnarray}
Then we recognize
\begin{eqnarray*}
&&\[B_\eps(x^\eps)\left(\frac{\partial}{\partial x_{i+1}}-\frac{\partial}{\partial x_{i}}\right) +
D_\eps(x^\eps)  \left ( \frac{\partial }{\partial x_{i}}-\frac{\partial }{\partial x_{i+1}}\right  )^2\]F(x^\eps)
\\
&=&
\left({\cal L}_{i,i+1}^{ABEP(\si,k)} F\right) (x^\epsi)
\end{eqnarray*}
which ends the proof of the lemma by the smoothness of $F$ and because by assumption, $x^\epsi\to x$.
\epr

\noindent
The weak asymmetry limit can also be performed on the $q$-TAZRP. This yields a totally asymmetric deterministic system as described in the following theorem.

\bt[Weak asymmetry limit of $q$-TAZRP]
Fix $T>0$.
Let $\{y^\epsi(t):0\leq T\}$ denote the
$q^\epsi$-TAZRP,  $q_\epsi:=1-\si\epsi$, with generator \eqref{genT11} and  initial condition
$y^\epsi(0)$. Assume that
\be\label{initialco1}
\lim_{\epsi\to  0} \epsi y^\epsi(0)= y\in {\mathcal X}_L
\ee
Then as $\epsi\to 0$, the process $\{y^\epsi(t):0\leq t\leq T\}$ converges weakly on
path space to the Totally Asymmetric Deterministic Energy Process,  TADEP$(\si)$ with generator
\be\label{gen-det}
({\cal L}^{TADEP}_{i,i+1}f)(z) = - \(\frac{1-e^{2\sigma z_{i+1}}}{2\sigma}\) \left(\frac{\partial}{\partial z_i} -\frac{\partial}{\partial z_{i+1}} \right)f(z), \qquad f: \mathbb{R}_+^L\to \mathbb R
\ee
initialized  from the configuration $y$.
\et
\bpr
The proof is analogous to the proof of Theorem \ref{teo1}
\epr
\subsection{Reversible measure of the ABEP$(\sigma,k)$}
\begin{theorem}[ABEP$(\sigma,k)$ reversible measures]
\label{basicproptheorem2}
\noindent
For all $L\in\mathbb{N}, L\geq 2$, the ABEP$(q,k)$ on ${\mathcal X}_L$ with closed boundary conditions
admits a family (labeled by $\gamma >-4\sigma k$) of reversible product measures with marginals given by
\be
\label{stat-meas-ABEP}
\mu_i(x_i):=\frac{1}{{\cal Z}^{(\gamma)}_i}  \, (1-e^{-2\sigma x_i})^{(2k-1)} e^{-(4\sigma k i+\gamma)x_i}\qquad\qquad x_i \in \mathbb R^+
\ee
for $i \in \Lambda_L$ and
\be \label{Zzz}
{\cal Z}_i^{{(\gamma)}} = \frac 1 {2\sigma}\; \text{Beta}\(2ki+\frac\gamma {2\si}, 2k\)
\ee
\end{theorem}
\bpr
The  adjoint of the generator of the ABEP$(\sigma,k)$ is given by
\begin{equation}
\({\cal L}^{ABEP^{(\sigma,k)}}_{(L)}\)^* =\sum_{i=1}^{{L-1}} \({\cal L}^{ABEP}_{i,i+1}\)^*
\end{equation}
with
\begin{eqnarray}
\({\cal L}^{ABEP}_{i,i+1}\)^*f&=&\frac{1}{4 \sigma^2} \, \(\frac{\partial}{\partial x_i}-\frac{\partial}{\partial x_{i+1}}\)^2 \,
\bigg(\(1-e^{-2\sigma x_i}\)\(e^{2\sigma x_{i+1}}-1\)f\bigg)\nonumber \\
&& \hskip-3cm -\frac{1}{2\sigma}  \(\frac{\partial}{\partial x_{i+1}}-\frac{\partial}{\partial x_{i}}\) \bigg(\left\{\(1-e^{-2\sigma x_i}\)\(e^{2\sigma x_{i+1}}-1\)+2k \[\(1-e^{-2\sigma x_i}\)-\(e^{2\sigma x_{i+1}}-1\)\]\right\}f\bigg) \nonumber
\end{eqnarray}
Let $\mu$ be a  product measure with $\mu(x)=\prod_{i=1}^L \mu_i(x_i)$, then in order for $\mu$ to be a stationary measure it is sufficient to impose that  the conditions
\begin{eqnarray}
&&\frac{1}{4 \sigma^2} \, \(\frac{\partial}{\partial x_{i+1}}-\frac{\partial}{\partial x_{i}}\) \,
\(1-e^{-2\sigma x_i}\)\(e^{2\sigma x_{i+1}}-1\)\mu(x)\nonumber \\
 &&-\, \frac{1}{2\sigma}  \, \left\{\(1-e^{-2\sigma x_i}\)\(e^{2\sigma x_{i+1}}-1\)+2k \[\(1-e^{-2\sigma x_i}\)-\(e^{2\sigma x_{i+1}}-1\)\]\right\} \mu(x) =0\nonumber
\end{eqnarray}
are satisfied for any $i \in \{1,\ldots, L-1\}$. This is true if and only if
 \begin{eqnarray}\label{cond}
 \frac{\mu_i'(x_i)}{\mu_i(x_i)}- 2 \sigma \, \frac{2k -e^{-2\sigma x_i}}{1-e^{-2\sigma x_i}}+\sigma =
  \frac{\mu_{i+1}'(x_{i+1})}{\mu_{i+1}(x_{i+1})}+ 2 \sigma \, \frac{e^{2\sigma x_{i+1}}-2k}{e^{2\sigma x_{i+1}}-1}-\sigma
 \end{eqnarray}
 for any $x_i, x_{i+1} \in \mathbb{R}^+$. The conditions \eqref{cond} are verified
 if and only if the marginals
 $\mu_i(x)$ are of the form   \eqref{stat-meas-ABEP} for some $\gamma \in \mathbb R$,  ${\cal Z}_i^{(\gamma)}$ is a normalization constant, and the constraint $\gamma >-4\si k$ is imposed in order to assure the integrability of $\mu(\cdot)$ on ${\mathcal X}_L$.  Thus we have proved that the product measure with marginal \eqref{stat-meas-ABEP} are stationary. One can also verify that for any $f:{\mathcal X}_L\to \mathbb R$
 $$
{\cal L}^{ABEP}f  =\frac 1 \mu \( {\cal L}^{ABEP}\)^*(\mu f)
 $$
 which then implies that the measure is reversible.
\epr
\br
In the limit $\si\to 0$ the reversible product measure of ABEP$(\sigma,k)$
converges to a product of Gamma distributions with shape parameter $2k$ and
scale parameter $1/\gamma$, which are the reversible homogeneous measures of
the BEP$(k)$ \cite{cggr}.
In the case $\si\neq 0$ the reversible product measure  of ABEP$(\sigma,k)$
has a decreasing  average profile (see Proposition \ref{decrease}).

\er
\subsection{Transforming the  ABEP$(\si,k)$ to BEP$(k)$}
In this subsection we show that the ABEP$(\sigma,k)$, which is an asymmetric process, can be mapped via a global change of coordinates
to the BEP(k) process which is symmetric. Here we focus on the analytical aspects of such $\sigma$-dependent mapping. In Section \ref{otto}
we will show that this map induces a conjugacy at the level of the underlying $\mathfrak{su}(1,1)$ algebra. This implies that the ABEP$(q,k)$ generator has a classical (i.e. non deformed) $\mathfrak{su}(1,1)$ symmetry. This is a remarkable because ABEP$(q,k)$ is a bulk-driven non-equilibrium process with non-zero average current (as it has been shown in Proposition \ref{non-zero})
and yet is generator is an element of the classical  $\mathfrak{su}(1,1)$ algebra.

\bd[Partial energy]\label{E}
We define the partial energy functions $E_i: {\mathcal X}_L\to \mathbb R_+$, $i \in \{1, \ldots, L+1\}$
\begin{equation}
 E_i(x):=\sum_{\ell = i}^L x_\ell, \qquad \text{for } i \in \Lambda_L \qquad \text{and } \quad E_{L+1}(x)=0.
\end{equation}
We also define the total energy $E: {\mathcal X}_L\to \mathbb R_+$ as 
$$ \qquad E(x):=E_1(x).$$
\ed
\bd[The mapping $g$]
\label{g-map}
We define the map $g:  {\mathcal X}_L \to {\mathcal X}_L$
\begin{equation}
\label{transfo}
g(x):=(g_i(x))_{i\in \Lambda_L} \qquad \text{with} \quad g_i(x):= \frac{e^{-2\sigma E_{i+1}(x)}-e^{-2\sigma E_i(x)}}{2\sigma}
\end{equation}
\ed
\noindent
Notice that $g$ does not have  full range, i.e.  $g[{\mathcal X}_L] \neq {\mathcal X}_L$. Indeed
\begin{equation}
E(g(x))= \frac 1 {2 \sigma} \(1-e^{-2\sigma E(x)}\)\le \frac1{2\sigma}
\end{equation}
so that in particular $g[{\mathcal X}_L]\subseteq \{x\in {\mathcal X}_L : E(x)\le 1/2\sigma\}$.
Moreover  $g$ is a bijection from $\mathcal{X}_L$
to $g[\mathcal{X}_L]$. Indeed, for  $z \in g[{\mathcal X}_L]$ we have
\begin{equation}\label{Invg}
(g^{-1}(z))_i=\frac{1}{2\si} \ln \left\{\frac{1-2\si \sum_{j=i+1}^L z_j}{1-2\si \sum_{j=i}^L z_j}\right\}
\end{equation}

\bt[Mapping from ABEP$(\sigma,k)$ to BEP$(k)$]
\label{map-theo}
Let $X(t)=(X_i(t))_{i \in \Lambda_L}$ be the ABEP$(\sigma, k)$ process starting from $X(0)=x$, then the process $Z(t):=(Z_i(t))_{i \in \Lambda_L}$ defined by the change of variable $Z(t):=g(X(t))$ is the BEP$(k)$ with initial condition $Z(0)=g(x)$.
\et

\begin{proof}
It is  sufficient to prove that, for any $f: {\mathcal X}_L\to \mathbb R_+$ smooth, $x \in  {\mathcal X}_L$  and $g$ defined above
\begin{equation}
\label{map-eq}
\[\mathcal L_{i,i+1}^{\text{BEP}} f\](g(x))= [\mathcal L_{i,i+1}^{\text{ABEP}} (f \circ g)](x)
\end{equation}
for any $i \in \Lambda_L$.
Define $F:=f \circ g$, then
\begin{eqnarray}
\label{yuh}
&&[\mathcal L^{\text{ABEP}} (f \circ g)](x) = [\mathcal L^{\text{ABEP}} (F)](x) =\\
& = & \frac 1{4\sigma^2}\, (1-e^{-2\sigma x_i})(e^{2\sigma x_{i+1}}-1)\(\frac{\partial}{\partial{x_{i+1}}}-\frac{\partial}{\partial{ x_{i}}}\)^2F(x)\nonumber \\
& +  & \frac{1}{2\sigma}\, \bigg\{ (1-e^{-2\sigma x_i})(e^{2\sigma x_{i+1}}-1)+2k\(2-e^{-2\sigma x_i} - e^{2\sigma x_{i+1}}\)\bigg\}\(\frac{\partial}{\partial{x_{i+1}}}-\frac{\partial}{\partial{x_{i}}}\)F(x)\nonumber
\end{eqnarray}
The computation of the Jacobian of $g$
\begin{equation}
\frac{\partial g_j}{\partial x_i}(x)=
\left\{
\begin{array}{ll}
-2\sigma g_j(x) & \text{for} \: j \le i-1\\
e^{-2\sigma E_j(x)}& \text{for} \: j = i\\
0 &\text{for} \: j \ge i+1\\
\end{array}
\right.
\end{equation}
implies that
\begin{equation}
\(\frac{\partial}{\partial x_{i+1}}-\frac{\partial}{\partial x_{i}}\)g_j(x)=
\left\{
\begin{array}{ll}
0 & \text{for} \: j \le i-1\\
-e^{-2\sigma E_{i+1}(x)}& \text{for} \: j = i\\
e^{-2\sigma E_{i+1}(x)}& \text{for} \: j = i+1\\
0 &\text{for} \: j \ge i+2\\
\end{array}
\right.
\end{equation}
and
\begin{equation}
\label{qaz}
\(\frac{\partial }{\partial x_{i+1}}-\frac{\partial }{\partial x_{i}}\)F(x)= e^{-2\sigma E_{i+1}(x)}\[\(\frac{\partial }{\partial z_{i+1}}-\frac{\partial }{\partial z_{i}}\)f\](g(x))
\end{equation}
\begin{eqnarray}
\label{wsx}
\(\frac{\partial }{\partial x_{i+1}}-\frac{\partial }{\partial x_{i}}\)^2F(x)&=& -2\sigma e^{-2\sigma E_{i+1}(x)}\[\(\frac{\partial }{\partial z_{i+1}}-\frac{\partial }{\partial z_i}\)f\](g(x)) \nn \\
\hskip-1cm&+&  e^{-4\sigma E_{i+1}(x)}\[\(\frac{\partial }{\partial z_{i+1}}-\frac{\partial }{\partial z_{i}}\)^2f\](g(x)).
\end{eqnarray}
Then, using \eqref{qaz} and \eqref{wsx}, \eqref{yuh} can be rewritten as
\begin{eqnarray}
&&[\mathcal L_{i,i+1}^{\text{ABEP}} (f \circ g)](x) =\nn \\
 &&\hskip1cm =  \frac 1{4\sigma^2}\, (1-e^{-2\sigma x_i})(e^{2\sigma x_{i+1}}-1)
 e^{-4\sigma E_{i+1}(x)} \[\(\frac{\partial }{\partial z_{i+1}}-\frac{\partial }{\partial z_{i}}\)^2f\](g(x))
\nonumber \\
 &&\hskip1cm+\bigg\{ 2\sigma  +  \frac{1}{2\sigma}\, \bigg( (1-e^{-2\sigma x_i})(e^{2\sigma x_{i+1}}-1)+2k\(2-e^{-2\sigma x_i} - e^{2\sigma x_{i+1}}\)\bigg)
\bigg\}e^{-2\sigma E_{i+1}(x)} \nn \\
\hskip1cm && \hskip10cm
\cdot \[\(\frac{\partial }{\partial z_{i+1}}-\frac{\partial }{\partial z_{i}}\)f\](g(x)) \nn
\end{eqnarray}
Simplifying, this gives
\begin{eqnarray}
&&[\mathcal L_{i,i+1}^{\text{ABEP}} (f \circ g)](x) =\nn \\
\hskip1cm &&\hskip1cm=  \bigg\{\frac{e^{-2\sigma E_{i+1}(x)}-e^{-2\sigma E_i(x)}}{2\sigma} \cdot \frac{e^{-2\sigma E_{i+2}(x)}-e^{-2\sigma E_{i+1}(x)}}{2\sigma}
\[\(\frac{\partial }{\partial z_{i+1}}-\frac{\partial }{\partial z_{i}}\)^2f\](g(x)) \nn \\
&& \hskip1cm- \frac{k}{\sigma} \(e^{-2\sigma E_{i}(x)}-2e^{-2\sigma E_{i+1}(x)}+e^{-2\sigma E_{i+2}(x)}\)
\[\(\frac{\partial }{\partial z_{i+1}}-\frac{\partial }{\partial z_{i}}\)f\](g(x)) \nn \\
\hskip1cm&&\hskip1cm =\[\mathcal L_{i,i+1}^{\text{BEP}} f\](g(x))\nn
\end{eqnarray}
\end{proof}

The ABEP$(\sigma,k)$ has a single conservation law given by the
total energy $E(x) = \sum_{i\in\Lambda_L} x_i$. As a consequence there exists an infinite
family of invariant measures which is hereafter described.
\bp[Microcanonical measure of ABEP$(\sigma,k)$]\label{CanMes}
The stationary measure of the ABEP$(\sigma,k)$ process on $\Lambda_L$
with given total energy $E$
is unique and is given by the inhomogeneous  product
measure with marginals \eqref{stat-meas-ABEP}
conditioned to a total energy $E(x)=E$.
More explicitly
\be
\label{mu-e}
d\mu^{(E)}(y) = \frac{\prod_{i=1}^L \mu_i(y_i) \mathbf{1}_{\{\sum_{i\in\Lambda_L} y_i = E\}}dy_i }{\int \ldots \int \prod_{i=1}^L \mu_i(y_i) \mathbf{1}_{\{\sum_{i\in\Lambda_L} y_i = E\}}dy_i}
\ee
\ep
\bpr
We start by observing that the stationary measure of the BEP$(k)$ process on $\Lambda_L$
with given total energy $\mathcal{E}$
is unique and is given by a  product of i.i.d.
Gamma  random variable $(X_i)_{i\in\Lambda_L}$
with shape parameter $2k$
conditioned to $\sum_{i\in\Lambda_l} X_i=\mathcal{E}$.
This is a consequence of duality between
BEP$(k)$ and SIP$(k)$ processes \cite{GKRV}.
Furthermore, an explicit computation shows that
the reversible measure of ABEP$(\sigma,k)$
conditioned to energy $E$ are transformed by the mapping $g$
(see Definition \ref{g-map})
to the stationary measure of the BEP$(k)$ with
energy $\mathcal{E}$ given by
$$\mathcal E= \frac{1}{2\sigma}(1-e^{-2\si E})\;.$$
The uniqueness for ABEP$(\sigma,k)$  follows
from the uniqueness for BEP$(\sigma,k)$
and the fact that $g$ is a bijection from $\mathcal{X}_L$
to $g[\mathcal{X}_L]$.
\epr

\subsection{The algebraic structure of ABEP$(\sigma,k)$}
\label{otto}

First we recall from \cite{GKRV} that the BEP$(k)$ generator can be written
in the form
\be
{\cal L}^{BEP(k)} = \sum_{i = 1}^{L-1} \left( K^{+}_i K_{i+1}^{-} + K^{-}_i K_{i+1}^{+} - K^{o}_i K_{i+1}^{o} + 2k^2 \right)
\ee
where
\begin{eqnarray}
\label{rep}
K^{+}_i & = & z_i \\
K^{-}_i & = & z_i \frac{\partial^2}{\partial z_i^2} + 2k \frac{\partial}{\partial z_i} \nonumber \\
K_i^{o} & = & z_i \frac{\partial}{\partial z_i} + k \nonumber
\end{eqnarray}
is a representation of the classical $\mathfrak{su}(1,1)$ algebra.
We show here that the ABEP$(\sigma,k)$ has the same algebraic structure.
This is proved by using a representation of $\mathfrak{su}(1,1)$
that is conjugated to \eqref{rep} and is given by
\be
\label{k-tilde}
\tilde{K}_i^{a}  = C_g \circ K_i^{a}  \circ  C_{g^{-1}}  \qquad\qquad \text{with} \quad a\in\{+,-,o\}
\ee
where $g$ is the function of Definition \ref{g-map} and
$$
(C_{g^{-1}} f)(x) = (f \circ g^{-1})(x)
$$
$$
(C_{g} f)(x) = (f \circ g)(x)\;.
$$
Explicitly one has
\be
(\tilde{K}_i^{a} f)(x) = (K_i^{a} f\circ g^{-1})(g(x)) \qquad\qquad \text{with} \quad a\in\{+,-,o\}
\ee

\bt[Algebraic structure of ABEP$(\sigma,k)$]
The generator of the ABEP$(\sigma,k)$ process is written as
\be
{\cal L}^{ABEP(\sigma, k)} = \sum_{i = 1}^{L-1} \left( \tilde{K}^{+}_i \tilde{K}_{i+1}^{-} + \tilde{K}^{-}_i \tilde{K}_{i+1}^{+} - \tilde{K}^{o}_i \tilde{K}_{i+1}^{o} + 2k^2 \right)
\ee
where the operators $\tilde{K}_i^{a}$ with $a\in\{+,-,o\}$ are defined in \eqref{k-tilde} and provide a representation
of the $\mathfrak{su}(1,1)$ Lie algebra.
\et
\bpr
The proof is a consequence of the following two results:
\be
\label{rewrite}
{\cal L}^{ABEP(\sigma,k)}  = C_g \circ {\cal L}^{BEP(k)}  \circ  C_{g^{-1}}
\ee
and the operators $\tilde{K}_i^{a}$ with $a\in\{+,-,o\}$ satisfy the commutation relations
of the $\mathfrak{su}(1,1)$ algebra. The first property is an immediate consequence
of Theorem \ref{map-theo}, as Eq. \eqref{rewrite} is simply a rewriting of Eq. \eqref{map-eq} by using
the definition of $C_g$ and $C_{g^{-1}}$. The second property can be obtained
by the following elementary Lemma, which implies that
the commutation relations
of the $\tilde{K}_i^{a}$ operators with $a\in\{+,-,o\}$ are the same of the
${K}_i^{a}$ operators with $a\in\{+,-,o\}$.
\epr
\bl
Consider an operator $A$ working on function $f:{\mathcal X}_L \to \mathbb{R}$ and let $g: {\mathcal X}_L \to X\subset {\mathcal X}_L$
be a bijection. Then defining
$$
\tilde{A} = C_g \circ A  \circ  C_{g^{-1}}
$$
we have that $A \to \tilde{A}$ is an algebra homomorphism.
\el
\bpr
We need to verify that
$$
\widetilde{A+B} = \tilde{A} + \tilde{B} \qquad\qquad \text{and} \qquad\qquad \widetilde{AB} = \tilde{A}\tilde{B}
$$
The first is trivial, the second is proved as follows
$$
\widetilde{AB} =  C_g \circ AB  \circ  C_{g^{-1}} = \left( C_g \circ A \circ C_{g^{-1}} \right) \circ \left( C_g \circ B  \circ  C_{g^{-1}}\right) = \tilde{A}\tilde{B}
$$
As a consequence
$$
\widetilde{[A,B]} = [\tilde{A},\tilde{B}]\;.
$$
\epr

\section{The Asymmetric KMP process, AKMP$(\sigma)$}
\label{akmp-section}
\subsection{Instantaneous Thermalizations}\label{instant}
The procedure of instantaneous thermalization has been introduced in \cite{GKRV}.
We consider a generator of the form
\be
\label{genera}
\mathcal L=\sum_{i} \mathcal L_{i,i+1}
\ee
where $\mathcal L_{i,i+1}$ is such that, for any initial
condition $(x_i,x_{i+1})$, the corresponding process converges to
a unique stationary distribution  $\mu_{(x_i,x_{i+1})}$.
\bd[Instantaneous thermalized process]
The instantaneous thermalization of the process with generator $\mathcal L$ in \eqref{genera}
is defined to be the process with generator
$$
\mathcal A=\sum_{i} \mathcal A_{i,i+1}
$$
where
\begin{eqnarray}
\label{A}
 \mathcal A_{i,i+1} f
 & = & \lim_{t\to\infty} (e^{t\mathcal L_{i,i+1}} f- f) \\
 & = & \int [f(x_1,\ldots,x_{i-1},y_i,y_{i+1},x_{i+2}, \ldots, x_L) - f(x_1,\ldots, x_L)] d\mu_{(x_i,x_{i+1})}(y_i,y_{i+1}) \nonumber
\end{eqnarray}
\ed
In words, in the process with generator $\mathcal A$ each edge $(i,i+1)$ is
updated at rate one, and after update its variables  are
replaced by a sample of the
stationary distribution of the process with generator $\mathcal L_{i,i+1}$
starting from $(x_i, x_{i+1})$.
Notice that, by definition, if a measure is stationary for the process with generator $\mathcal L_{i,i+1}$ then
it is also stationary for the  process with generator $\mathcal A_{i,i+1}$.\\

\noindent
An  example of thermalized processes is the Th-BEP$(k)$ process, where the local redistribution
rule is
\be
\label{betaaa}
(x, y) \to (B(x+y), (1-B)(x+y))
\ee
with $B$ a Beta$(2k,2k)$ distributed random variable \cite{CGGR}.
In particular for $k=1/2$ this gives the KMP process \cite{KMP} that has a uniform redistribution rule on $[0,1]$.
Among discrete models we mention the Th-SIP$(k)$ process where the redistribution rule is
\be
(n,m) \to (R, n+m -R)
\ee
where $R$ is Beta-Binomial$(n+m,2k,2k)$. For $k=1/2$ this
corresponds to discrete uniform distributions on $\{0,1,\ldots, n+m\}$.
Other examples are described in \cite{CGGR}.
In the following we introduce the asymmetric version
of these redistribution models.

\subsection{Thermalized Asymmetric Inclusion process Th-ASIP$(q,k)$}
The instantaneous thermalization limit of the Asymmetric Inclusion process  is obtained as follows.
Imagine on each bond $(i,i+1)$ to run the ASIP$(q,k)$ dynamics for an infinite amount of time.
Then the total number of particles   on the bond will be redistributed according to the stationary
measure on that bond, conditioned to conservation of the total number of particles of the bond.
We consider the independent random variables  $(M_1,\ldots,M_L)$ distributed according
to the stationary measure of the ASIP$(q,k)$ at  equilibrium.
Thus  $M_i$ and $M_{i+1}$ are distributed according to
\be
\label{stat-meas1}
p^{(\alpha)}_i(\eta_i):=\mathbb{P^{(\alpha)}}(M_i = \eta_i) = \frac{\alpha^{\eta_i}}{Z^{(\alpha)}_i} \,{\binom{\eta_i+2k-1}{\eta_i}_q} \cdot  q^{4ki\eta_i}  \qquad\qquad \eta_i \in \mathbb N
\ee
and
\be
\label{stat-meas2}
p^{(\alpha)}_{i+1}(\eta_{i+1}):=\mathbb{P^{(\alpha)}}(M_{i+1} = \eta_{i+1}) = \frac{\alpha^{\eta_{i+1}}}{Z^{(\alpha)}_{i+1}} \,{\binom{\eta_{i+1}+2k-1}{\eta_{i+1}}_q} \cdot  q^{4k(i+1)\eta_{i+1}}  \qquad\qquad \eta_{i+1} \in \mathbb N
\ee
for some $\alpha\in [0,q^{-(2k+1)})$.
\noindent
Hence the distribution of $M_i$, given that the sum is fixed to $M_i + M_{i+1} =n+m$ has  the following  probability
mass function:
\begin{eqnarray}
\label{qbbb}
\nu^{ASIP}_{q,k}(r\, | \, n+m)&:=&\mathbb{P}(M_i=r\: |\:M_i+M_{i+1} =n+m)  \\
&=&  \frac{p^{(\alpha)}_i(r) p^{(\alpha)}_{i+1}(n+m-r)}{\sum_{l=0}^{n+m} p^{(\alpha)}_i(l) p^{(\alpha)}_{i+1}(n+m-l)} \nn \\
&=& \widetilde{\mathcal C}_{q,k}(n+m)\, q^{-4kr}\, \binom{r+2k-1}{r}_q \cdot \binom{2k+n+m-r-1}{n+m-r}_q \nn
\end{eqnarray}
where $r\in \mathbb N$ and $ \widetilde{\mathcal C}_{q,k}(n+m)$ is a normalization constant.
\bd[Th-ASIP$(q,k)$ process]
The Th-ASIP$(q,k)$ process on $\Lambda_L$ is defined as the thermalized discrete
process with state space $\Omega_L$ and local redistribution rule
\be
(n,m) \to (R_q, n+m -R_q)
\ee
where $R_q$ has a $q$-deformed Beta-Binomial$(n+m,2k,2k)$  distribution
with mass function \eqref{qbbb}.
The generator of this process is given by
\begin{eqnarray}\label{genTHasip}
&&{\cal L}^{ASIP(q,k)}_{th}f(\eta) \nonumber
\\
&= &
\sum_{i=1}^{L-1} \sum_{r=0}^{\eta_i+\eta_{i+1}}\[f(\eta_1, \dots, \eta_{i-1},r, \eta_i+\eta_{i+1}-r, \eta_{i+2}, \dots, \eta_L)-f(\eta)\]\;\nu^{ASIP}_{q,k}(r\, | \,\eta_i+\eta_{i+1})
\nonumber\\
\end{eqnarray}
\ed
\subsection{Thermalized Asymmetric Brownian energy process  Th-ABEP$(\sigma,k)$.}
We define the instantaneous thermalization limit of the Asymmetric Brownian Energy process  as follows.
On each bond we run the ABEP$(\si,k)$ for an infinite time.
Then the energies on the bond will be redistributed according to the stationary
measure on that bond, conditioned to the conservation of the total energy of the bond.
If we take two independent random variables  $X_i$ and $X_{i+1}$ with  distributions as in \eqref{stat-meas-ABEP}, i.e.
\be
\label{stat-meas-ABEP1}
\mu_i(x_i):=\frac{1}{{\cal Z}^{(\gamma)}_i}  \, (1-e^{-2\sigma x_i})^{(2k-1)} e^{-(4\sigma k i+\gamma)x_i}\qquad\qquad x_i \in \mathbb R^+
\ee
\be
\label{stat-meas-ABEP2}
\mu_{i+1}(x_{i+1}):=\frac{1}{{\cal Z}^{(\gamma)}_{i+1}}  \, (1-e^{-2\sigma x_{i+1}})^{(2k-1)} e^{-(4\sigma k(i+1)+\gamma)x_{i+1}}\qquad\qquad x_{i+1} \in \mathbb R^+
\ee
then the distribution of  $X_i$, given  the sum  fixed to $X_i + X_{i+1} =E$, has
density
\begin{eqnarray}
p(x_i|X_i+X_{i+1} =E) &=&\frac{\mu_i(x_i)\mu_{i+1}(E-x_{i})}{\int_0^E \mu_i(x)\mu_{i+1}(E-x) \, dx} \nonumber \\
&=&
\mathcal C_{\sigma,k}(E) \, e^{4\sigma k x_i} \[\(1-e^{-2\sigma x_i}\)\(1-e^{-2\sigma (E-x_i)}\)\]^{2k-1} \nn
\end{eqnarray}
where $\mathcal C_{\sigma,k}(E) $ is a normalization constant. Equivalently, let $W_i := X_i/E$, then $W_i$ is a random
variable taking values on $[0,1]$. Conditioned to $X_i+X_{i+1}=E$, its  density is given by
\begin{eqnarray}
\label{plutooo}
\nu_{\sigma,k}(w|E)= \widehat{\mathcal C}_{\sigma,k}(E)\; e^{2\sigma E w}
\left\{\(e^{2\sigma E w}-1\)\(1-e^{-2\sigma E(1-w)}\)\right\}^{2k-1}
\end{eqnarray}
with
\begin{eqnarray}
\quad  \widehat{\mathcal C}_{\sigma,k}(E):= \int_0^1 e^{2\sigma E w}
\left\{\(e^{2\sigma E w}-1\)\(1-e^{-2\sigma E(1-w)}\)\right\}^{2k-1} \; dw
\end{eqnarray}
\bd[Thermalized ABEP$(\sigma,k)$]\label{DEF}
The Th-ABEP$(\sigma,k)$ process on $\Lambda_L$ is defined as the thermalized
process with state space ${\mathcal X}_L$ and local redistribution rule
\be
(x,y) \to (B_{\sigma}(x+y), (1-B_{\sigma})(x+y) )
\ee
where $B_{\sigma}$ has a  distribution with density function $\nu_{\sigma,k}(\cdot|x+y)$
in \eqref{plutooo}. Thus the
generator of Th-ABEP$(\si,k)$  is given by
\begin{eqnarray}\label{ThABEP}
&& {\cal L}^{ABEP(\sigma,k)}_{th}f(x) =
 \\ 
&& \hskip.1cm=
\sum_{i=1}^{L-1} \int_{0}^1 \[f(x_1,\ldots, w(x_i+x_{i+1}), (1-w)(x_i+x_{i+1}),\ldots,x_L) - f(x)\]\, \nu_{\sigma,k}(w|x_i+x_{i+1}) \,dw \nn
\end{eqnarray}
\ed

\noindent
In the limit $\sigma\to0$, the conditional density $\nu_{0^+,k}(\cdot|E)$ does not depend on $E$, and for any $E\ge 0$ we recover the Beta$(2k,2k)$ distribution
with density
\begin{eqnarray}
\nu_{0^+,k}(w|E)= \frac{1}{\text{Beta}(2k,2k)}\[w(1-w)\]^{2k-1}\;.
\end{eqnarray}
Then the generator ${\cal L}^{ABEP(0^+,k)}_{th}$ coincides with the generator of the thermalized Brownian Energy process Th-BEP$(k)$ defined in equation (5.13) of \cite{cggr}.\\

\noindent
The redistribution rule with the random variable $B_{\sigma}$ in Definition \ref{DEF} is truly asymmetric,
meaning that - on average -  the energy is moved to the left.
\bp
\label{decrease}
Let $B_\si$ be the random variable  on $[0,1]$ distributed with density \eqref{plutooo}, then $\mathbb E[B_\si]\ge\frac{1}{2}$. As a consequence $B_\si$ and $1-B_\si$ are not equal in distribution and for $(X_1,\ldots,X_L)$ distributed according to the reversible product measure $\mu$ of  ABEP$(\si,k)$ defined in \eqref{stat-meas-ABEP}, we have that the energy profile is decreasing, i.e.
\be
\label{general}
\mathbb E_\mu[X_i]\ge\mathbb E_\mu[X_{i+1}], \qquad \forall \, i \in \{1, \ldots, L-1\}\;.
\ee
\ep
\bpr
Let $X = (X_1,X_2)$ be a two-dimensional random vector taking values in $\mathcal{X}_2$ distributed according to the microcanonical  measure $\mu^{(E)}$ of ABEP$(\si,k)$
with fixed total energy $E\ge 0$, defined in \eqref{mu-e}. Then, from Definition \ref{DEF},
\begin{equation}
(X_1,X_2) \,{\buildrel d \over =}\, (EB_\sigma, E(1-B_\sigma)) \qquad \text{with} \qquad B_\sigma \sim  \nu_{\sigma,k}(\cdot|E)
\end{equation}
Then, as already remarked in the proof of Proposition \ref{CanMes}, $Z:=g(X)$ with $g(\cdot)$  as in Definition \ref{g-map}  is a two-dimensional random variable taking values in $g[\mathcal{X}_2]\subset \mathcal{X}_2$
and distributed according to the microcanonical  measure of BEP$(k)$ with fixed total energy $\mathcal E= \frac{1}{2\sigma}(1-e^{-2\si E})$.
It follows from \eqref{betaaa} that
\begin{equation}
g(X)  \,{\buildrel d \over =}\, (\mathcal EB, \mathcal E(1-B)) \qquad \text{with} \qquad B \sim \text{Beta}(2k,2k)\;.
\end{equation}
Then, by \eqref{Invg} we have
\begin{equation}
(1-B_\sigma)E =  (g^{-1}(Z))_2= \frac{1}{2\sigma} \ln \left\{\frac{1}{1-2\si (1-B)\mathcal E}\right\}
\end{equation}
and therefore 
\be
B_{\sigma} = {1+}\frac{1}{2\sigma E} \ln \left(1-B (1-e^{{-}2\sigma E})\right)
\ee
Put $2\si E=1$ without loss of generality, for simplicity.
Then to prove that $\mathbb{E}[B_\sigma] > 1/2$ we have to prove that
$$
\E (1+ \ln(1-B(1-e^{-1})))\geq \frac12
$$
Defining  $a=1-e^{-1}$ we then have to prove that
\be
\E(-\ln(1-aB))\leq \frac12
\ee
It is useful to write
$$
-\ln(1-aB)= \sum_{n=1}^\infty \frac{a^n B^n}{n}
$$
and remark that for a $Beta(\alpha,\alpha)$ distributed
$B$ one has
$$
\E(B^n)= \prod_{r=0}^{n-1} \frac{\alpha+r}{2\alpha+r}\;.
$$
So we have to prove that
$$
\psi(\alpha, a):=\sum_{n=1}^\infty \frac{a^n}{n}\prod_{r=0}^{n-1} \frac{\alpha+r}{2\alpha+r}<1/2
$$
First consider {the limit $\alpha\to\infty$} then we find
$$
{\lim_{\alpha \to \infty}\phi(\alpha,a)}=\sum_{n=1}^\infty \frac{a^n}{2^n n}= -\ln \left(1-\frac12 (1-e^{-1})\right)=-\ln \left(\frac12 +\frac{e^{-1}}{2}\right)\approx 0.379<1/2
$$
Next remark
when $\alpha=0$ the $B$ is distributed like  $\frac12 \delta_0+\frac12 \delta_1$ which gives
$$
\E(-\ln(1-aB))=-\frac12 \ln (e^{-1})=\frac12
$$
Now we prove that $\psi$ is monotonically decreasing in $\alpha$. 
To see this notice that
$$
\frac{d}{d\alpha} \frac{\alpha+r}{2\alpha+r} =\frac{-r}{(2\alpha+r)^2} <0
$$
So the derivative
$$
\frac{d}{d\alpha}\psi(\alpha, a)= \sum_{n=1}^\infty\sum_{r'=0}^{n-1}\frac{a^n}{n} \left(\prod_{r=0, r\not= r'}^{n-1}\frac{\alpha+r}{2\alpha+r}\right)\frac{-r'}{(2\alpha+r)^2}<0
$$
Therefore $\psi(\alpha, a)$ is monotonically decreasing in $\alpha$ and 
$\psi(\alpha, a)\leq \frac12$. Thus the claim $\mathbb{E}[B_\sigma] > 1/2$ is proved.

%
Now let $X=(X_1,X_2)$ be a two-dimensional r.v. distributed according to the profile measure $\mu$ defined in \eqref{stat-meas-ABEP} with $L=2$
and with abuse of notation let $\nu_{\sigma,k} \left[B_\si|E\right] = \mathbb E \left[B_\si\right]$. Then we can write $X=(E \, B_\sigma,  E(1-B_\sigma)) $ where now $E$ is a  random
variable. We have
\begin{eqnarray}
\mathbb E_\mu \left[X_2\right]&=& \mathbb E_\mu \left[\mathbb E_\mu\left[X_2| \;  E \right] \right]= \mathbb E_\mu \left[\mathbb E_\mu\left[E(1-B_\sigma)| \;  E \right] \right]=\mathbb E_\mu \left[E \; \nu_{\sigma,k} \left[(1-B_\si)|\;  E\right] \right] \nn \\
&\le & \mathbb E_\mu \left[E \; \nu_{\sigma,k} \left[B_\si|\;  E\right] \right]= \mathbb E_\mu \left[\mathbb E_\mu\left[X_1| \;  E \right] \right]=\mathbb E_\mu \left[X_1\right]
\end{eqnarray}
The proof can be easily generalized to the case $L\ge 2$, yielding \eqref{general}.
\epr
\vskip.2cm
\noindent
For $k=1/2$ and $\si\to 0$ the Th-ABEP$(\si,k)$ is exactly the KMP process
 \cite{KMP}.  For $k=1/2$ and $\si>0$
 \begin{eqnarray}
\nu_{\sigma,1/2}(w|E)= \frac{2\sigma E}{e^{2\sigma E}-1} e^{2\sigma E w}, \qquad w\in [0,1]
\end{eqnarray}
The Th-ABEP$(\si,\frac12)$ can therefore
be considered as the natural asymmetric analogue of the KMP process.
This justifies the following definition.
\bd[AKMP$(\sigma)$ process]
We define the Asymmetric KMP with asymmetry parameter $\sigma\in\mathbb{R}_+$
on $\Lambda_L$ as the process with generator given by:
\begin{eqnarray}
&&{\cal L}^{AKMP(\sigma)}f(x)= \sum_{i=1}^{L-1} \left\{ \frac{2\sigma (x_i+x_{i+1})}{e^{2\sigma (x_i+x_{i+1})}-1}\;\cdot \right.
 \nonumber \\
&&
\left. \cdot \; \int_{0}^1  \[f(x_1,\ldots, w(x_i+x_{i+1}), (1-w)(x_i+x_{i+1}),\ldots,x_L) - f(x)\]\,
 e^{2\sigma w(x_i+x_{i+1})} \,dw \right\}\nn
\end{eqnarray}
\ed


\section{Duality relations}\label{dualsect}
In this section we derive various duality properties of the processes
introduced in the previous sections. We start by recalling the definition of {\em duality}.
\bd
\label{standard}
Let  $\{X_t\}_{t\ge0}$,  $\{\widehat{X}_t\}_{t\ge0}$ be two  Markov processes with state spaces  $\Omega$ and $\widehat{\Omega}$ and $D: \Omega\times \widehat{\Omega}\to\R$ a bounded
measurable function. The processes  $\{X_t\}_{t\ge0}$,  $\{\widehat{X}_t\}_{t\ge0}$ are said to be dual with respect to $D$  if
\be\label{standarddualityrelation1}
\E_x \big[D(X_t, \widehat{x})\big]=\widehat{\mathbb{E}}_{\widehat{x}} \big[D(x, \widehat{X}_t)\big]\;
\ee
for all $x\in\Omega, \widehat{x}\in \hat{\Omega}$ and $t>0$. In (\ref{standarddualityrelation1})
$\E_x $ is the expectation with respect to the
 law of the $\{X_t\}_{t\ge0}$ process started at $x$, while $\widehat{\mathbb{E}}_{\widehat{x}} $
denotes  expectation with respect to the law of the $\{\widehat{X}_t\}_{t\ge0}$ process
initialized at $\widehat{x}$.
\ed

\subsection{Self-duality of ASIP$(q,k)$}
\label{subsect51}
The basic duality relation is the self-duality of  ASIP$(q,k)$.
This self-duality property is derived from a symmetry of the underlying Hamiltonian
which is a sum of co-products of the Casimir operator. In \cite{CGRS} this construction
was achieved for the algebra ${\mathcal{U}}_q(\mathfrak{su}(2))$, and  from the Hamiltonian a Markov
generator was constructed via a positive ground state. Here the construction
and consequent symmetries is analogous, but for the algebra ${\mathcal{U}}_q(\mathfrak{su}(1,1))$.
For the proof of the following Theorem we refer to
Section \ref{D}, where we implement the
steps of \cite{CGRS}
for the algebra ${\mathcal{U}}_q(\mathfrak{su}(1,1))$.

\bt[Self-duality of the finite ASIP$(q,k)$]
\label{mainself}
The ASIP$(q,k)$ on $\Lambda_L$ with closed boundary conditions
is self-dual with the following self-duality function
\be
\label{dualll}
D_{(L)}( \eta, \xi )=
\prod_{i =1}^L  \frac{\binom{\eta_i}{\xi_i}_q}{\binom{\xi_i+2k-1}{\xi_i}_q}\,
\cdot \,
q^{(\eta_i-\xi_i)\left[2\sum_{m=1}^{i-1}\xi_m +\xi_i\right]- 4 k i \xi_i}
\cdot \mathbf 1_{\xi_i \le \eta_i}
\ee
or, equivalently,\be
\label{duallll}
D_{(L)}( \eta, \xi )=
\prod_{i =1}^L  \frac{(q^{2(\eta_i-\xi_i+1)}; q^2)_{\xi_i}}{(q^{4k}; q^2)_{\xi_i}} \cdot q^{(\xi_i-4ki+2N_{i+1}(\eta))\xi_i}
\cdot \mathbf 1_{\xi_i \le \eta_i}
\ee
with $(a;q)_m$ as defined in \eqref{Poch} and
\begin{equation}
\label{current}
 N_i(\eta):= \sum_{k= i}^L \eta_{{k}}\;.
\end{equation}
\et

\br\label{ConfXi}
For $n\in\mathbb{N}$, let $\xi^{(\ell_1, \ldots, \ell_n)}$ be the configurations with
$n$ particles located at sites $\ell_1, \ldots, \ell_n$.
Then for the configuration $\xi^{(\ell)}$ with one particle at site $\ell$
\begin{equation}
\label{one-dual}
D(\eta,\xi^{(\ell)})  = \frac{q^{-(4k\ell+1)}}{q^{2k}-q^{-2k}} \, \cdot (q^{2 N_\ell(\eta)}-q^{2 N_{\ell+1}(\eta)})
\end{equation}
and, more generally, for the configuration $\xi^{(\ell_1, \ldots, \ell_n)}$  with
$n$ particles at  sites $\ell_1, \ldots, \ell_n$ with $\ell_i \neq \ell_j$
$$
D(\eta,\xi^{(\ell_1, \ldots, \ell_n)}) = \frac{q^{-4k\sum_{m=1}^n \ell_m-n^2}}{(q^{2k}-q^{-2k})^n} \, \cdot \prod_{m=1}^n(q^{2 N_{\ell_m}(\eta)}-q^{2 N_{\ell_m+1}(\eta)})
$$
\er

\vskip.2cm
\noindent
The duality relation with duality function \eqref{duallll} makes sense in the limit $L \to \infty$. Indeed, if $N_i(\eta)=\infty$ for some $i$, then $\lim_{L \to \infty}D_{(L)}(\eta,\xi)=0$ for all $\xi$ with $\xi_i \neq 0$. If the initial configuration $\eta\in \Omega_\infty$ has a finite number of particles at the right of the origin, then from the duality relation, we deduce that  it remains like this for all later times $t>0$, which implies that $N_\ell(\eta_t)<\infty$ for all $t\ge 0$. Conversely, if $\eta$ is such that $N_0(\eta)=\infty$, then  $N_0(\eta_t)=\infty$ for all later times because, from the duality relation, $\mathbb E_\xi \left[D(\eta, \xi_t)\right]=0$ for all $t>0$. To extract some non-trivial informations from the duality relation in the infinite volume case, a suitable  renormalization is required (see Section \ref{Infi}).

\subsection{Duality between  ABEP$(\si,k)$ and  SIP$(k)$}
\label{subsect52}
We remind the reader that in the limit of zero asymmetry $q\to 1$
the ASIP$(q,k)$ converges to the SIP$(k)$.
Therefore from the self-duality of
ASIP$(q,k)$,
and the fact that the ABEP$(\si,k)$ arises
as a limit of ASIP$(q,k)$ with $q\to 1$, a duality between ABEP$(\si,k)$
and SIP$(k)$ follows.
\bt[Duality ABEP$(\sigma,k)$ and SIP$(k)$]
\label{main}
The ABEP$(\sigma,k)$ on $\Lambda_L$ with closed boundary conditions
is dual to the SIP$(k)$ on $\Lambda_L$ with closed boundary conditions, with the following self-duality function
\be
\label{dualABEP}
D_{(L)}^\sigma(x,\xi)=\prod_{i \in \Lambda_L}\frac{\Gamma(2k)}{\Gamma(2k+\xi_i)}\, \left(\frac{e^{-2\sigma E_{i+1}(x)}- e^{-2\sigma E_{i}(x)}}{2\sigma}\right)^{\xi_i}
\ee
with $E_i(\cdot)$ the partial energy function defined in Definition \ref{E}.
\et
\bpr
The duality function in \eqref{dualABEP} is related to the duality function between BEP$(k)$ and SIP$(k)$, $D^0_{(L)}(x, \eta)$ (see e.g. Section 4.1 of \cite{cggr}) by the following relation
\begin{equation}
D_{(L)}^\sigma(x,\xi)= D^0_{(L)}(g(x), \eta)
\end{equation}
where  $g(\cdot)$ is the map defined in \eqref{g-map}. Thus, omitting the subscript $(L)$ in the following, from \eqref{rewrite} we have
\begin{eqnarray}
\left[\mathcal{L}^{\text{ABEP}(\si,k)}D^\sigma(\cdot,\eta)\right](x)&=&\left[\mathcal{L}^{\text{ABEP}(\si,k)}\(D^0(\cdot,\eta)\circ g\)\right](x)\nn \\
& = & \left[\mathcal{L}^{\text{BEP}(k)}D^0(\cdot,\eta)\right](g(x)) \nn \\
&=&\left[\mathcal{L}^{\text{SIP}(k)}D^0(g(x),\cdot)\right](\eta)\nn \\& =&\left[\mathcal{L}^{\text{SIP}(k)}D^\sigma(x,\cdot)\right](\eta)
\end{eqnarray}
and this proves the Theorem.
\epr
\br
In the limit as $\sigma\to 0$ one recovers the duality $D^0_{(L)}(\cdot,\cdot)$ between  BEP$(k)$ and  SIP$(k)$.
However it is remarkable here that for finite $\sigma$ there is duality between a bulk driven
asymmetric process, the ABEP$(\sigma,k)$,
and an equilibrium symmetric process, the SIP$(k)$. Indeed, the asymmetry is hidden in the
duality function. This is somewhat reminiscent of the dualities between systems
with reservoirs and absorbing systems \cite{cggr}, where also the source of non-equilibrium,
namely the different parameters of the reservoirs has been moved to the duality
function.
\er

\vskip.3cm
\noindent
The following proposition explains how $D_{(L)}^\sigma(x,\xi)$ arises as the limit of ASIP$(q,k)$ self-duality function for $q=1-N^{-1} \sigma$, $N \to \infty$.
\bp For any fixed $L\ge 2$ we have
\begin{equation}\label{limit}
 \lim_{N \to \infty}\(\frac{\sigma}{N}\)^{|\xi|} \, D_{(L)}^{\text{ASIP}(1-\si/N,k)}(\lfloor Nx \rfloor,\xi)=  D_{(L)}^{\text{ABEP}(\si,k)}(x,\xi)
\end{equation}
where $D_{(L)}^{\text{ASIP}(q,k)}(\eta,\xi)$ denotes the self-duality function of ASIP$(q,k)$ defined in \eqref{duallll} and $D_{(L)}^{\text{ABEP}(\si,k)}(x,\xi)$ denotes the duality function defined in \eqref{dualABEP}.
\ep

\bpr
 Let
\be
N:=|\eta|:=\sum_{i=1}^L\eta_i, \qquad q=1-\frac{\sigma}{N}, \qquad  x:= N^{-1}\eta,
\ee
then
\be
 D_{(L)}^{\text{ASIP}(q,k)}(\eta, \xi )=
\prod_{i =1}^L  \frac{[\eta_i]_q[\eta_i-1]_q \ldots[\eta_i-\xi_i+1]_q}{[2k+\xi_i-1]_q[2k+\xi_i-2]_q\ldots[2k]_q}\,
\cdot \,
q^{(\eta_i-\xi_i)\left[2\sum_{m=1}^{i-1}\xi_m +\xi_i\right]- 4 k i \xi_i}
\cdot \mathbf 1_{\xi_i \le \eta_i}
\ee
Now, for any $m$
\begin{eqnarray}
[\eta_i-m]_{1-\frac{\sigma}{N}}&=&[Nx_i - m]_{1-\frac{\sigma}{N}}\nonumber \\
&=& \frac{N}{2\sigma}\[e^{\sigma x_i}-e^{-\sigma x_i}+ O(N^{-1})\] \nonumber \\
&=& \frac{N}{\sigma}\sinh(\sigma x_i)+ O(1)
\end{eqnarray}
hence
\be
\prod_{m=0}^{\xi_i-1}[Nx_i-m]_{1-\frac{\sigma}{N}}= \(\frac{N}{\sigma}\, \sinh (\sigma x_i)+  O(1)\)^{\xi_i}
\ee
On the other hand
\be
[2k+m]_{1-\frac{\sigma}{N}}=2k+m+  O(N^{-1}) \qquad \text{thus} \qquad \prod_{m=0}^{\xi_i-1} [2k+m]_{1-\frac{\sigma}{N}}= \frac{\Gamma(2k+\xi_i)}{\Gamma(2k)}+ O(N^{-1})
\ee
finally, let $f_i(\xi):=2\sum_{m=1}^{i-1}\xi_m +\xi_i$ and $g_i(\xi):=-\xi_i\[2\sum_{m=1}^{i-1}\xi_m +\xi_i\]-4ki\xi_i$ we have
\be
q^{\eta_i f_i(\xi)}= \(1-\frac{\sigma}{N}\)^{Nx_i f_i(\xi)}=e^{-\sigma x_i f_i(\xi)}+  O(N^{-1}), \quad \text{and} \quad
q^{g(\xi)}= \(1-\frac{\sigma}{N}\)^{g(\xi)}= 1 +  O(N^{-1})
\ee
then
 \eqref{limit} immediately follows.
\epr


\subsection{Duality for the instantaneous thermalizations}
\label{subsect53}
{In this section we will prove that the self-duality of ASIP$(q,k)$ and the duality between ABEP$(\si,k)$ and SIP$(k)$ imply duality properties also for the thermalized models.}
\bp
If a process $\{\eta(t): t\geq0\}$ with generator $\mathcal L=\sum_{i=1}^{L-1} \mathcal L_{i,i+1}$ is dual
to a process $\{\xi(t) : t\geq 0\}$ with generator  $\widehat{\mathcal L} = \sum_{i=1}^{L-1} \widehat{\mathcal L}_{i,i+1}$
with duality function $D(\cdot, \cdot)$ in such a way
that for all $i$
\begin{equation*}
\left[\mathcal L_{i,i+1} D(\cdot, \xi)\right](\eta)=[\widehat{\mathcal L}_{i,i+1} D(\eta, \cdot)](\xi)
\end{equation*}
then, if  the instantaneous thermalization processes of $\eta_t$, resp. $\xi_t$ both exist,  they are each other's
dual with the same duality function $D(\cdot, \cdot)$.
\ep
\bpr
Let $\mathcal A$, resp. $\widehat {\mathcal A}$ be  the generators of  the instantaneous thermalization of $\eta_t$, resp. $\xi_t$, then, from \eqref{A} we know that
\begin{equation*}
\mathcal A= \sum_{i \in \Lambda_L} \mathcal A_{i,i+1}, \qquad \mathcal A_{i,i+1}= \lim_{t\to\infty} (e^{t \mathcal L_{i,i+1}}-I)
\end{equation*}
and
\begin{equation*}
\widehat{\mathcal A}= \sum_{i \in \Lambda_L} \widehat{\mathcal A}_{i,i+1}, \qquad \widehat{\mathcal A}_{i,i+1}= \lim_{t\to\infty} (e^{t \widehat{\mathcal L}_{i,i+1}}-I)
\end{equation*}
where $I$ denotes identity and where the exponential $e^{t\mathcal L_{i,i+1}}$ is the semigroup generated by $\mathcal L_{i,i+1}$ in the sense of the Hille Yosida theorem.
Hence we immediately
obtain that
\begin{equation*}
\Big[(e^{t \mathcal L_{i,i+1}}-I)D(\cdot, \xi)\Big](\eta)= \Big[(e^{t\widehat{\mathcal L}_{i,i+1}}-I) D(\eta, \cdot)\Big](\xi)
\end{equation*}
which proves the result.
\epr

\vskip.2cm
\noindent
As a consequence of this Proposition we obtain duality between the thermalized
ABEP$(q,k)$ and the thermalized SIP$(k)$ as well as self-duality of
the thermalized ASIP$(q,k)$.
\bt\label{thermaldual}
\vskip.2cm
\noindent
\begin{itemize}
\item[a)] The Th-ASIP$(q,k)$ with generator \eqref{genTHasip} is self-dual with self-duality
function given by \eqref{dualll}.
\item[b)]
The Th-ABEP$(\sigma,k)$ with generator \eqref{ThABEP}  is dual, with duality function \eqref{dualABEP} to the Th-SIP$(k)$ in $\Lambda_L$ whose generator is given by
\begin{eqnarray}\label{ThSIP}
&&{\cal L}^{SIP(k)}_{th}f(\xi)=  \\
& & = \sum_{i=1}^{L-1} \sum_{r=0}^{\xi_i+\xi_{i+1}}\[f(\xi_1, \dots, \xi_{i-1},r, \xi_i+\xi_{i+1}-r, \xi_{i+2}, \dots, \xi_{L})-f(\xi)\]\;\nu^{SIP}_{k}(r\, | \,\xi_i+\xi_{i+1}) \nn
\end{eqnarray}
where $\nu^{SIP}_{k}(r\, | \,n+m)$ is the probability density of a Beta-Binomial distribution of parameters $(n+m,2k,2k)$.
\end{itemize}
\et

\vskip.2cm

\br
For $k=1/2$ \eqref{ThSIP} gives
the KMP-dual, i.e., the asymmetric KMP  has the same dual as
the symmetric KMP, but of course with different $\si$-dependent
duality function given by
\be
\label{dualAKMP}
D_{(L)}^{\text{AKMP}(\si)}(x,\xi)=\prod_{i \in \Lambda_L}\frac{1}{\xi_i!}\, \left(\frac{e^{-2\sigma E_{i+1}(x)}- e^{-2\sigma E_{i}(x)}}{2\sigma}\right)^{\xi_i}
\ee

\er

\section{Applications to exponential moments of currents}\label{computations}

The  definition of the ASIP$(q,k)$ process on the infinite lattice requires
extra conditions on the initial data. Indeed, when the total number of particles
is infinite,  there is the possibility of the appearance of singularities,
since a single site can accommodate an unbounded number of
particles. By self-duality we can however  make sense of expectations
of duality functions in the infinite volume limit. This is the aim of the next section.

\subsection{Infinite volume limit for ASIP$(q,k)$}\label{Infi}
In this section we approximate an infinite-volume configuration
by a finite-volume configuration and we appropriately renormalize
the self-duality function to avoid divergence in the thermodynamical
limit.
\bd[Good infinite-volume configuration]\label{Good}
\vskip.1cm
\noindent
\begin{itemize}
\item[a)] We say that $\eta \in \mathbb{N}^{\mathbb{Z}}$ is a {\em ``good infinite-volume configuration''} for ASIP$(q,k)$ iff for  $\eta^{(L)}\in \mathbb{N}^{\mathbb{Z}}$, $L \in \mathbb{N}$,  the restriction of $\eta$ on $[-L,L]$, i.e.
\begin{equation}\label{EtaL}
\eta^{(L)}_i=
\left\{
\begin{array}{ll}
\eta_i & \text{for} \quad  i \in [-L,L]\\
0 & \text{otherwise}
\end{array}
\right.
\end{equation}
the limit
\begin{equation}
\lim_{L \to \infty} \prod_{i \in \mathbb Z}q^{-2\xi_i N_{i+1}(\eta^{(L)})}\;\mathbb E_\xi \left[D(\eta^{(L)},\xi(t))\right]
\end{equation}
exists and is finite for all $t \ge 0$ and for any $\xi \in \mathbb{N}^{\mathbb{Z}}$ finite (i.e. such that $\sum_{i\in \mathbb Z}\xi_i<\infty$).
\item[b)] Let $\mu$ be a probability measure on $\mathbb{N}^{\mathbb{Z}}$, then we say that it is a {\em  ``good infinite-volume  measure''} for ASIP$(q,k)$ iff it concentrates on good infinite-volume configurations.
\end{itemize}
\ed
\bp\label{InfVol}
\vskip.1cm
\noindent
\begin{itemize}
\item[1)] If $\eta \in \mathbb{N}^{\mathbb{Z}}$ is a { ``good infinite-volume configuration''} for ASIP$(q,k)$ and  $\xi^{(\ell_1, \ldots, \ell_n)}$ is the configurations with
$n$ particles located at sites $\ell_1, \ldots, \ell_n\in \mathbb Z$, then the limit
\be
 \lim_{L\to\infty} \prod_{m=1}^n q^{-2N_{\ell_m +1}(\eta^{(L)})} \; \mathbb{E}_{\eta^{(L)}}\left[D(\eta(t),\xi^{(\ell_1, \ldots, \ell_n)})\right]
\ee
is well-defined for all $t\ge 0$ and is equal to
\begin{equation}
\lim_{L \to \infty} \prod_{m=1}^n q^{-2N_{\ell_m +1}(\eta^{(L)})} \;\mathbb E_{\xi^{(\ell_1, \ldots, \ell_n)}} \left[D(\eta^{(L)},\xi(t))\right]
\end{equation}
\item[2)] If $\eta \in \mathbb{N}^{\mathbb{Z}}$ is bounded, i.e.  $\sup_{i\in\mathbb{Z}} \eta_i <\infty$, then it is a ``good infinite-volume configuration''.
\item[3)] Let us denote by $\mathcal N_\lambda(t)$ a Poisson process of rate $\lambda>0$, and by $\mathbf E[\cdot]$ the expectation w.r. to its probability law. If $\mu$ is a probability measure on $\mathbb{N}^{\mathbb{Z}}$ such that for any $\lambda>0$ the expectation
\begin{equation}\label{mu}
\mathbb{E}_\mu \left[ \mathbf E \left[ e^{\sum_{i=1}^{\mathcal N_\lambda(t)}\eta_{\ell +i}}\right]\right]
\end{equation}
is finite for all $t \ge 0$ and for any $\ell \in \mathbb Z$, then $\mu$ is a ``good infinite-volume  measure''.
\end{itemize}
\ep

\bpr
\begin{itemize}
\item[1)]
 If $\eta \in \mathbb{N}^{\mathbb{Z}}$ is a good infinite volume configuration, then the duality relation with duality function \eqref{duallll} makes sense after the following renormalization:
\begin{equation}\label{QUI}
 \mathbb{E}_{\eta^{(L)}}\left[D(\eta(t),\xi^{(\ell_1, \ldots, \ell_n)})\right] \prod_{m=1}^n q^{-2N_{\ell_m +1}(\eta^{(L)})}=
 \mathbb{E}_{\xi^{(\ell_1, \ldots, \ell_n)}}\left[D(\eta^{(L)},\xi(t))\right] \prod_{m=1}^n q^{-2N_{\ell_m +1}(\eta^{(L)})}
 \end{equation}
 then the first statement of the Theorem follows after taking the limit as $L \to \infty $ of \eqref{QUI}.
 \item[2)] Let $\xi$ be a  finite configuration in $\mathbb{N}^{\mathbb{Z}}$.
 We prove that for any bounded  $\eta \in \mathbb{N}^{\mathbb{Z}}$  the family of functions
  \begin{eqnarray}
  \mathcal S_L(t):=\prod_{i \in \mathbb Z}q^{-2\xi_i N_{i+1}(\eta^{(L)})}\;\mathbb E_\xi \left[D(\eta^{(L)},\xi(t))\right], \qquad L \in \mathbb N \label{Start}
  \end{eqnarray}
  is uniformly bounded.
 Without loss of generality we can suppose that $\xi=\xi^{(\ell_1, \ldots, \ell_n)}$, for some $\{\ell_1, \ldots, \ell_n\}\subset \mathbb Z$, $n \in \mathbb N$. Moreover we  denote by $(\ell_1(t), \ldots, \ell_n(t))$ the positions of the $n$ ASIP$(q,k)$  walkers starting at time $t=0$ from $(\ell_1, \ldots, \ell_n)$. We then have $\xi(t)=\xi^{(\ell_1(t), \ldots, \ell_n(t))}$, and   \begin{eqnarray}
\mathcal S_L(t) &=&  \prod_{m=1}^n q^{-2N_{\ell_m +1}(\eta^{(L)})} \; \mathbb{E}_{\xi^{(\ell_1, \ldots, \ell_n)}}\left[D(\eta^{(L)},\xi(t))\right]= \nn \\
 && =\mathbb{E}_{\xi^{(\ell_1, \ldots, \ell_n)}}\Bigg[
\prod_{i =1}^L  \frac{(q^{2(\eta^{(L)}_i-\xi_i(t)+1)}; q^2)_{\xi_i(t)}}{(q^{4k}; q^2)_{\xi_i(t)}} \cdot q^{\xi_i^2(t)}
\cdot \mathbf 1_{\xi_i(t) \le \eta^{(L)}_i} \cdot \nn \\
 &&\cdot
  \prod_{m =1}^n   q^{-4k\ell_m(t)+2[N_{\ell_m(t)+1}(\eta^{(L)})-N_{\ell_m+1}(\eta^{(L)})]} \Bigg]\;.\nn
  \end{eqnarray}
As a consequence, since
\begin{equation}
 (q^{2(\eta-\xi+1)}; q^2)_{\xi}  \cdot q^{\xi^2}
\cdot \mathbf 1_{\xi\le \eta} \le 1
\end{equation}
and
\be
\sup_{\ell \le n}\frac{1}{(q^{4k}; q^2)_{\xi}}\le c
\ee
for some $c>0$,
 we have that there exists $C>0$ such that
\begin{equation}\label{CCC}
\big|\mathcal S_L(t) \big|\le C \; \mathbb{E}_{\xi^{(\ell_1, \ldots, \ell_n)}}\left[ \prod_{m =1}^n  q^{-4k\ell_m(t)+2[N_{\ell_m(t)+1}(\eta^{(L)})-N_{\ell_m+1}(\eta^{(L)})]}\right]
\end{equation}
for all $L\in \mathbb N$, $t\ge 0$. Then, from  the Cauchy-Schwarz inequality, in order to find an upper bound  for \eqref{CCC}, it is sufficient to find an upper bound  for
\begin{eqnarray}
  s_{L,m}(t):=\mathbb{E}_{\xi^{(\ell_1, \ldots, \ell_n)}}\left[q^{\kappa\{-4k\ell_m(t)+2[N_{\ell_m(t)+1}(\eta^{(L)})-N_{\ell_m+1}(\eta^{(L)})]\}}\right]\nn
\end{eqnarray}
for any fixed $m \in \{1, \ldots, n\}$ and $\kappa \in \mathbb N$. Now, let  $M:=\sup_{i\in\mathbb{Z}} \eta_i <\infty$, then
$$
\big|N_{\ell_m(t)+1}(\eta^{(L)})-N_{\ell_m+1}(\eta^{(L)})\big|\le M |\ell_m(t)-\ell_m|
$$
hence  there exists $C',\omega>0$  such that
\begin{eqnarray}\label{QUO}
\big|s_{L,m}(t) \big|\le C' \; \mathbb{E}_{\xi^{(\ell_1, \ldots, \ell_n)}}\left[e^{\omega|\ell_m(t)-\ell_m|}\right]
\end{eqnarray}
for any $L\in \mathbb N$, $t \ge 0$. Since $\xi(t)$ has a finite number of particles, for each $m \in \{1, \ldots, n\}$ the process $|\ell_m(t)-\ell_m|$ is stochastically dominated by a Poisson process $\mathcal N(t)$ with parameter
\be \label{la}
\lambda:= \max_{0\le\eta,\eta'\le n}\{q^{\eta-\eta'+(2k-1)} [\eta]_q [2k+\eta']_q\} \vee \max_{0\le\eta,\eta'\le n}\{  q^{\eta-\eta'-(2k-1)} [2k+\eta]_q [\eta']_q
\}
\ee
then the right hand side of \eqref{QUO} is less or equal than
\begin{equation}
\mathbf E\[e^{\omega \mathcal N(t)}\]= e^{-\lambda t}\sum_{i=0}^{\infty} e^{\omega i} \frac{(\lambda t)^i}{i!}<\infty\;.
\end{equation}
This proves that $\mathcal S_L(t)$ is uniformly bounded.
 \item[3)] Suppose that the probability measure $\mu$ satisfies \eqref{mu}. Then, in order to prove that it is a ``good'' measure, it is sufficient to show that
 \be
\lim_{L \to \infty} \mathbb E_\mu \left[ \prod_{i \in \mathbb Z}q^{-2\xi_i N_{i+1}(\eta^{(L)})}\;\mathbb E_\xi \left[D(\eta^{(L)},\xi(t))\right]\right]<\infty \label{Start1}
\ee
By exploiting  the same arguments  used in  the proof of item 2), we claim that, in order to prove \eqref{Start1} it is sufficient to show that for each fixed $m=1, \ldots, n$, $\kappa>0$, the function
\begin{eqnarray}
  \Theta_{L,m}(t):=\mathbb E_{\mu}\left[\mathbb{E}_{\xi^{(\ell_1, \ldots, \ell_n)}}\left[q^{\kappa\{-4k\ell_m(t)+2[N_{\ell_m(t)+1}(\eta^{(L)})-N_{\ell_m+1}(\eta^{(L)})]\}}\right]\right]
\end{eqnarray}
is uniformly bounded.
We have that
\begin{eqnarray}
 && \Theta_{L,m}(t)=\nn \\
 &&=\mathbb E_{\mu}\left[\mathbb{E}_{\xi^{(\ell_1, \ldots, \ell_n)}}\left[q^{-4\kappa k\ell_m(t)}
  \left(q^{-2\kappa \sum_{i= \ell_m+1}^{\ell_m(t)} \eta_i^{(L)}}\mathbf 1_{\ell_m<\ell_m(t)}+q^{2\kappa \sum_{i= \ell_m(t)+1}^{\ell_m} \eta_i^{(L)} }\mathbf 1_{\ell_m(t)<\ell_m}\right)\right]\right]\nn \\
  &&\le  \mathbb E_{\mu}\left[\mathbb{E}_{\xi^{(\ell_1, \ldots, \ell_n)}}\left[q^{-4\kappa k\ell_m(t)}
  \left(q^{-2\kappa \sum_{i= 1}^{\ell_m(t)-\ell_m} \eta_{i+\ell_m}^{(L)}}\mathbf 1_{\ell_m<\ell_m(t)}+1\right)\right]\right]\;.\nn \end{eqnarray}
Then the result follows  as in proof of item 2) from the fact that the process $\ell_m(t)-\ell_m$ is stochastically dominated by a Poisson process of rate $\lambda$ \eqref{la}, and from the hypothesis \eqref{mu}.
\end{itemize}
\epr

\noindent
Later on, if we write expectations in the infinite volume we always refer to the
limiting procedure described above. Namely, for a ``good infinite-volume configuration'' $\eta \in \mathbb{N}^{\mathbb{Z}}$, with an abuse of notation we will write
\be\label{notation}
 \prod_{m=1}^n q^{-2N_{\ell_m +1}(\eta)}\;\mathbb{E}_{\eta}\left[D(\eta(t),\xi^{(\ell_1, \ldots, \ell_n)})\right]
:= \lim_{L\to\infty} \prod_{m=1}^n q^{-2N_{\ell_m +1}(\eta^{(L)})}
\; \mathbb{E}_{\eta^{(L)}}\left[D(\eta(t),\xi^{(\ell_1, \ldots, \ell_n)})\right]
\ee
and
\begin{equation}\label{notation1}
 \prod_{m=1}^n q^{-2N_{\ell_m +1}(\eta)} \;\mathbb E_{\xi^{(\ell_1, \ldots, \ell_n)}} \left[D(\eta,\xi(t))\right]:=\lim_{L \to \infty} \prod_{m=1}^n q^{-2N_{\ell_m +1}(\eta^{(L)})} \;\mathbb E_{\xi^{(\ell_1, \ldots, \ell_n)}} \left[D(\eta^{(L)},\xi(t))\right]
\end{equation}

\subsection{$q$-exponential moment of the current of ASIP$(q,k)$}
We start by defining the current for the ASIP$(q,k)$ process on $\mathbb{Z}$.
\bd[Current] Let $\{\eta(t), \; t \ge 0\}$ be a c\`adl\`ag trajectory on the infinite-volume  configuration space $\mathbb N^{\mathbb Z}$, then
the total integrated current $J_i(t)$ in the time interval $[0,t]$ is defined as the net number of particles crossing the bond $(i-1,i)$
in the right direction. Namely, let $(t_i)_{i\in\mathbb{N}}$ be the sequence of the process jump times. Then
\be
J_i(t) = \sum_{k: t_k \in [0,t]}  \(\mathbf{1}_{\{\eta(t_k) = \eta(t_k^-)^{i-1,i}\}} - \mathbf{1}_{\{\eta(t_k) = \eta(t_k^-)^{i,i-1}\}}\)
\ee
\ed
\bl[Current]
The total integrated current of a c\`adl\`ag trajectory $(\eta(s))_{0\le s\le t}$ with $\eta(0)=\eta$ is given by
\be \label{J}
J_i(t)= N_i(\eta(t))-N_i(\eta):= \lim_{L \to \infty}\(N_i(\eta^{(L)}(t))-N_i(\eta^{(L)})\)
\ee
where $N_i(\eta)$ is defined in \eqref{current} and $\eta^{(L)}$ is defined in \eqref{EtaL}.
Moreover
\be\label{J0}
\lim_{i\to -\infty} J_i(t) = 0
\ee
\el
\bpr
\eqref{J} immediately follows from the definition of $J_i(t)$, whereas \eqref{J0} follows from the conservation of the total number of particles.
\epr

\bp[Current q-exponential moment via a dual walker]\label{Lemma:N}
Let $\eta \in \mathbb N^{\mathbb Z}$ a good infinite-volume configuration in the sense of Definition \ref{Good}, then the first $q$-exponential moment of the current when the process is started from
$\eta$ at time $t=0$ is given by
\be
\mathbb{E}_{\eta}\[q^{2J_i(t)}\] = q^{2(N(\eta)-N_i(\eta))} - \sum_{n=-\infty}^{i-1} q^{4kn} \; \mathbf E_n \left[ q^{-4km(t)}\(1-q^{-2\eta_{m(t)}}\)\, q^{2(N_{m(t)}(\eta)-N_i(\eta))}\right]
\label{Quii}
\ee
where
$m(t)$ denotes a continuous time asymmetric random walker on $\mathbb Z$
jumping left at rate $q^{-2k}[2k]_q $ and jumping right at rate $q^{2k}[2k]_q $
and $\mathbf E_i$ denotes the expectation
with respect to the law of $m(t)$ started at site $i\in\mathbb{Z}$ at time $t=0$.
Furthermore $N(\eta)-N_i(\eta) = \sum_{n<i} \eta_n$ and the first term
on the right hand side of \eqref{Quii} is zero when there are infinitely many particles
to the left of $i\in\mathbb{Z}$ in the configuration $\eta$.
\ep

\bpr
To prove \eqref{Quii}
we consider    the configuration $\xi^{(i)}\in \mathbb{N}^{\mathbb Z}$  with a single dual particle at site $i$. Since the ASIP$(q,k)$ is self-dual the dynamics of the single dual
particle is given an asymmetric random walk $m(t)$ on $\mathbb Z$  whose
rates are computed from the process definition and coincides
with those in the statement of the Proposition.
From \eqref{notation}, \eqref{notation1} and item 1) of Proposition \ref{InfVol} we have that
$$
q^{-2N_i(\eta)}\;\mathbb{E}_{\eta}\left[D(\eta(t),\xi^{(i)})\right] = \frac{q^{-(4ki+1)}}{q^{2k}-q^{-2k}} \; q^{-2N_i(\eta)}\,\mathbb{E}_{\eta}\left[q^{2 N_i(\eta(t))}-q^{2 N_{i+1}(\eta(t))}\right]
$$
is equal to
$$
q^{-2N_i(\eta)}\;\mathbb{E}_{\xi^{(i)}}\left[D(\eta,\xi^{(m(t))})\right]= q^{-2N_i(\eta)}\frac{q^{-1}}{q^{2k}-q^{-2k}} \; \mathbf{E}_{i}\left[q^{-4km(t)}(q^{2 N_{m(t)}(\eta)}-q^{2N_{m(t)+1}(\eta)})\right]
$$

\noindent
Then from \eqref{J} we get
\begin{eqnarray}\label{R1}
\mathbb{E}_{\eta}\left[q^{2 J_i(t)}\]
& = &
q^{-2\eta_i}\;\mathbb{E}_{\eta}\left[q^{2 J_{i+1}(t)}\right] \nonumber \\
& + &
q^{4ki}\; \mathbf{E}_{i}\left[q^{-4km(t)}(q^{2 (N_{m(t)}(\eta)-N_i(\eta))}-q^{2(N_{m(t)+1}(\eta)-N_i(\eta))})\right]
\end{eqnarray}
{By iterating the relation in \eqref{R1}, for any $n\ge 0$ we get
\begin{eqnarray}\label{R11}
&&\mathbb{E}_{\eta}\left[q^{2 J_{i+1}(t)}\]
= 
q^{2(N_{i-n}(\eta)-N_{i+1}(\eta))}\;\mathbb{E}_{\eta}\left[q^{2 J_{i-n}(t)}\right] +\nonumber \\
& &- 
\sum_{j=0}^n q^{2(N_{i-j}(\eta)-N_{i+1}(\eta))}q ^{4k(i-j)}\; \mathbf{E}_{i-j}\left[q^{-4km(t)}(q^{2 (N_{m(t)}(\eta)-N_{i-j}(\eta))}-q^{2(N_{m(t)+1}(\eta)-N_{i-j}(\eta))})\right]\;.\nn\\
\end{eqnarray}
By taking the limit $n\to \infty$ we get
\begin{eqnarray}
\mathbb{E}_{\eta}\left[q^{2 J_{i+1}(t)}\]&=&
\lim_{n\to \infty}q^{2(N_{i-n}(\eta)-N_{i+1}(\eta))}\;\mathbb{E}_{\eta}\left[q^{2 J_{i-n}(t)}\right] + \nn
\\
&-&
\sum_{j=0}^\infty q^{-2N_{i+1}(\eta)}q ^{4k(i-j)}\; \mathbf{E}_{i-j}\left[q^{-4km(t)}(q^{2 N_{m(t)}(\eta)}-q^{2N_{m(t)+1}(\eta)})\right]\nn
\end{eqnarray}
and using \eqref{J0} we obtain \eqref{Quii}.}
\epr

%
%
We continue with a lemma that is useful in the following.
\bl \label{Lemma:LDP}
Let $x(t)$ be the random walk on $\mathbb Z$ jumping to the right with rate $a\ge 0$ and to the left with rate $b\ge 0$, let $\alpha \in \mathbb R$, and $A \subseteq \mathbb R$ then
\be \label{LDP}
\lim_{t \to \infty} \frac 1 t \log \mathbf E_0 \[\alpha^{x(t)}\, \Big| \; x(t)\in A \]= \sup_{x \in A}\left\{x \log \alpha -{\cal I}(x)\right\} - \inf_{x \in A}{\cal I}(x)
\ee
with
\be \label{I}
{\cal I}(x)=(a+b)-\sqrt{x^2+4ab} + x \ln \(\frac{x+\sqrt{x^2+4ab}}{2a}\)
\ee
\el
\bpr
From large deviations theory \cite{denholla} we know that $x(t)/t$, conditioned to
$x(t)/t \in A$, satisfies a large
deviation principle with rate function
${\cal I}(x) - \inf_{x\in A} {\cal I}(x)$ where ${\cal I}(x)$ is given by
\be
{\cal I}(x):= \sup_{z}\left\{zx- \Lambda(z)\right\}
\ee
with
\be
\Lambda(z):= \lim_{t \to \infty} \frac 1 t \log \mathbb E\[e^{zx(t)}\]= a\(e^z-1\)+b\(e^{-z}-1\)
\ee
from which it easily follows \eqref{I}. The application
of Varadhan's lemma yields \eqref{LDP}.
\epr

\vskip.2cm
\br \label{Rem:LDP}
Let $m(t)$ be the random walk defined in Proposition \ref{Lemma:N}, then
\eqref{LDP} holds with
\be \label{I1}
{\cal I}(x)=[4k]_q-\sqrt{x^2+(2[2k]_q)^2} + x \log \left\{\frac{1}{2[2k]_qq^{2k}}\[x+\sqrt{x^2+(2[2k]_q)^2} \]\right\}
\ee
\er

\vskip.3cm

\noindent
We denote by  ${\mathbb E}^{\otimes \mu}$ the expectation of the ASIP$(q,k)$ process on $\mathbb{Z}$
initialized with the  homogeneous product measure on $\mathbb{N}^{\mathbb Z}$ with marginals $\mu$
at time 0, i.e.
$$
{\mathbb E}^{\otimes \mu}[f(\eta(t))]= \sum_\eta \left(\otimes_{i\in\mathbb{Z}} \mu(\eta_i)\right) \mathbb E_\eta[f(\eta(t))]\;.
$$

\bp[$q$-moment for product initial condition]
\label{current-prod}
Consider an homogeneous product probability measure $\mu$ on $\mathbb N$. Then,
 for the infinite volume ASIP$(q,k)$, we have
\be \label{Exp-mom-prod}
\mathbb{E}^{\otimes\mu}\[q^{2J_i(t)}\]= \mathbf E_0 \[\(\frac{q^{-4k}}{\lambda_q}\)^{m(t)} \mathbf 1_{m(t)\le 0}\]+
\mathbf E_0 \[q^{-4km(t)} \(\lambda_{1/q}^{m(t)}-\lambda_{1/q}+ \lambda_q^{-1}\)\mathbf 1_{m(t)\ge 1}\]
\ee
where $\lambda_y:= \sum_{n=0}^{\infty} y^n \mu(n)$
and $m(t)$ is the random walk defined in Proposition \ref{Lemma:N}.
In particular we have
\be\label{M}
\lim_{t \to\infty} \frac{1}{t} \log\mathbb E^{\otimes\mu}[q^{2J_i(t)}]= \sup_{x \ge 0}\left\{x \log M_q -{\cal I}(x)\right\} -
 \inf_{x \ge 0}{\cal I}(x)
\ee
with $M_q:=  q^{-4k}\lambda_{1/q}$ and ${\cal I}(x)$ given by \eqref{I1}.
\ep
\bpr
It is easy to check that an homogeneous product measure $\mu$ verifies the condition \eqref{mu} in Proposition \ref{InfVol}, thus it is a good infinite-volume probability measure in the sense of Definition \ref{Good}. For this reason we can apply Proposition  \ref{Lemma:N}, and  from \eqref{Quii} we have
\begin{eqnarray}
&&\mathbb{E}^{\otimes\mu}\[q^{2J_i(t)}\]= \int \otimes\mu(d\eta)\, \mathbb E_\eta\[q^{2J_i(t)}\]
\nn \\
&& = \int \otimes\mu(d\eta) q^{2(N(\eta)-N_i(\eta))}
+ \sum_{n=-\infty}^{i-1} q^{4kn} \; \int\otimes\mu(d\eta) \mathbf E_n \left[ q^{-4km(t)}\(q^{-2\eta_{m(t)}}-1\)\, q^{2(N_{m(t)}(\eta)-N_i(\eta))}\right] \nn
\end{eqnarray}
Since
\be
\int \otimes\mu(d \eta) q^{2(N_m(\eta)-N_i(\eta))} = \lambda_q^{i-m} \; \mathbf 1_{\{m \le i\}} + \lambda_{1/q}^{m-i} \; \mathbf 1_{\{m>i\}}
\ee
then, in particular, $\int \otimes\mu(d\eta) q^{2(N(\eta)-N_i(\eta))} = 0$ since $\lambda_q<1$,
where we recall the interpretation of $N(\eta)-N_i(\eta)$ from  Proposition \ref{Lemma:N}.
Hence
\begin{eqnarray}\label{Ct}
\mathbb{E}^{\otimes\mu}\[q^{2J_i(t)}\]&=&  \sum_{n=-\infty}^{i-1} q^{4kn} \;\sum_{m \in \mathbb Z} \mathbf P_n\(m(t)=m\) \, q^{-4km} \int \otimes\mu(d \eta) \left[q^{2(N_{m+1}(\eta)-N_i(\eta))}-q^{2(N_{m}(\eta)-N_i(\eta))}\right] \nn \\
&=& \(\lambda_q^{-1}-1\) A(t) + \(\lambda_{1/q}-1\)B(t)
\end{eqnarray}
with
\be
A(t):=\sum_{n\le i-1} q^{4kn} \sum_{m \le i} \mathbf P_n\(m(t)=m\) q^{-4km} \lambda_q^{i-m}
\ee
and
\be
B(t):=\sum_{n\le i-1} q^{4kn} \sum_{m \ge i+1} \mathbf P_n\(m(t)=m\) q^{-4km} \lambda_{1/q}^{m-i}
\ee
Now, let $\alpha:= q^{-4k}\lambda_q^{-1}$, then
\begin{eqnarray}
A(t)&=& \sum_{n\le i-1} q^{4kn} \lambda_q^{i}\sum_{m \le i} \mathbf P_n\(m(t)=m\) \alpha^m \nn \\
&=& \sum_{j\ge 1} \lambda_q^j \sum_{\bar{m} \le j} \mathbf P_0 \(m(t)=\bar{m}\) \alpha^{\bar{m}} \nn \\
&=& \sum_{\bar m \le 0} \alpha^{\bar m} \mathbf P_0 \(m(t)=\bar m\)  \sum_{j\ge 1} \lambda_q^j + \sum_{\bar m\ge 1} \alpha^{\bar m} \mathbf P_0\(m(t)=\bar m\)\sum_{j \ge \bar m} \lambda_q^j \nn \\
&=& \frac{1}{1-\lambda_q} \left\{ \lambda_q \, \mathbf E_0 \[\alpha^{m(t)}\, \mathbf 1_{m(t)\le 0}\]+ \mathbf E_0 \[q^{-4km(t)}\, \mathbf 1_{m(t)\ge 1}\]\right\} \label{At}
\end{eqnarray}
Analogously one can prove that
\be
B(t)=\frac{1}{\lambda_{1/q}-1}\left\{\mathbf E_0 \[\beta^{m(t)}\, \mathbf 1_{m(t)\ge 2}\]-\lambda_{1/q}\mathbf E_0 \[q^{-4km(t)}\, \mathbf 1_{m(t)\ge 2}\]\right\}\label{Bt}
\ee
with $\beta=q^{-4k}\lambda_{1/q}$ then \eqref{Exp-mom-prod} follows by combining \eqref{Ct}, \eqref{At} and \eqref{Bt}.
\vskip.2cm
\noindent
In order to prove \eqref{M} we use the fact that $m(t)$ has a Skellam distribution with parameters $([2k]_q q^{2k}t, [2k]_q q^{-2k}t)$,
i.e. $m(t)$ is the difference of two independent Poisson random variables with those parameters.  This
implies that
\be
 \mathbf E_0 \[\(\frac{q^{-4k}}{\lambda_q}\)^{m(t)} \mathbf 1_{m(t)\le 0}\]=  \mathbf E_0 \[\lambda_q^{m(t)} \mathbf 1_{m(t)\ge 0}\]\;. \nn
\ee
Then we can rewrite \eqref{Exp-mom-prod}
as
\begin{eqnarray} \label{Exp-mom-prod1}
\mathbb{E}^{\otimes\mu}\[q^{2J_i(t)}\]
& = &
\mathbf E_0 \[\lambda_q^{m(t)} \, \mathbf 1_{m(t)\ge 1}\] + \mathbf P_0\(m(t)=0\) \nn \\
&+ & \(\lambda_q^{-1}-\lambda_{1/q}\) \mathbf E_0 \[q^{-4km(t)} \mathbf 1_{m(t)\ge 1}\]+ \mathbf E_0 \[M_q^{m(t)}\mathbf 1_{m(t)\ge 1}\]\nn \\
& = & \mathbf E_0 \[M_q^{m(t)} \mathbf 1_{m(t)\ge 0}\]  \left(1+\mathcal E_1(t)+\mathcal E_2(t)+\mathcal E_3(t) +\mathcal E_4(t)\right)
\end{eqnarray}
with
\be
\mathcal E_1(t) :=\frac{\mathbf E_0 \[M_q^{m(t)}\mathbf 1_{m(t)\ge 1}\]}
{\mathbf E_0 \[M_q^{m(t)} \mathbf 1_{m(t)\ge 0}\]}, \qquad \mathcal E_2(t):=\frac{\mathbf P_0\(m(t)=0\)}{\mathbf E_0 \[M_q^{m(t)} \mathbf 1_{m(t)\ge 0}\]}
\nn
\ee
and
\be
\mathcal E_3(t):=\frac{\mathbf E_0 \[\lambda_q^{m(t)} \mathbf 1_{m(t)\ge 1}\]}{\mathbf E_0 \[M_q^{m(t)} \mathbf 1_{m(t)\ge 0}\]}, \qquad \mathcal E_4(t):=\frac{\(\lambda_q^{-1}-\lambda_{1/q}\)\mathbf E_0 \[q^{-4km(t)} \mathbf 1_{m(t)\ge 1}\]}{\mathbf E_0 \[M_q^{m(t)} \mathbf 1_{m(t)\ge 0}\]}
\ee
To identify the leading term in \eqref{Exp-mom-prod1} it remains to prove that, for each $i=1,2,3$ there exists $c_i>0$ such that
\be\label{Ee}
\sup_{t\ge 0}|\mathcal  E_i(t)| \le c_i
\ee
This would imply, making use of  Lemma \ref{Lemma:LDP}, the result in \eqref{M}. The bound in \eqref{Ee} is immediate for $i=1,2,3$. To prove it for $i=4$ it is sufficient to show that there exists $c>0$ such that
\be
\lambda_q^{-1}\mathbf E_0 \[q^{-4km(t)} \mathbf 1_{m(t)\ge 1}\] \le c \, \mathbf E_0 \[\(q^{-4k}\lambda_{1/q}\)^{m(t)}\mathbf 1_{m(t)\ge 1}\]\;.
\ee
This follows since there exists $m_*\ge 1$ such that for any $m \ge m_*$ $\lambda_q^{-1}\le \lambda_{1/q}^m$ and then
\begin{eqnarray}
\lambda_q^{-1}\mathbf E_0 \[q^{-4km(t)} \mathbf 1_{m(t)\ge 1}\]\le \lambda_q^{-1}
\mathbf E_0 \[q^{-4km(t)} \mathbf 1_{1\le m(t)< m_*}\]+ \mathbf E_0 \[q^{-4km(t)} \lambda_{1/q}^{m(t)} \mathbf 1_{m(t)\ge m_*}\]\nn \\
\le \lambda_q^{-1}
\mathbf E_0 \[q^{-4km(t)} \mathbf 1_{1\le m(t)}\]+ \mathbf E_0 \[q^{-4km(t)} \lambda_{1/q}^{m(t)} \mathbf 1_{m(t)\ge 1}\]\nn \\
\le \(1+\lambda_q^{-1}\)\mathbf E_0 \[\(q^{-4k} \lambda_{1/q}\)^{m(t)} \mathbf 1_{m(t)\ge 1}\]\;.
\end{eqnarray}
This concludes the proof.
\epr

\subsection{Infinite volume limit for ABEP$(\sigma,k)$}

\bd[Good infinite-volume configuration]\label{Good1}
\vskip.1cm
\noindent
\begin{itemize}
\item[a)] We say that $x \in \mathbb{R}_+^{\mathbb{Z}}$ is a {\em ``good infinite-volume configuration''} for ABEP$(\si,k)$ iff for  $x^{(L)}\in \mathbb{R}_+^{\mathbb{Z}}$, $L \in \mathbb{N}$,  the restriction of $x$ to $[-L,L]$, i.e.
\begin{equation}\label{xL}
x^{(L)}_i=
\left\{
\begin{array}{ll}
x_i & \text{for} \quad  i \in [-L,L]\\
0 & \text{otherwise}
\end{array}
\right.
\end{equation}
the limit
\begin{equation}
\lim_{L \to \infty} \prod_{i \in \mathbb Z}e^{2\si\xi_i E_{i+1}(x^{(L)})}\;\mathbb E_\xi \left[D^\si(x^{(L)},\xi(t))\right]
\end{equation}
exists and is finite for all $t \ge 0$ and for any $\xi \in \mathbb{N}^{\mathbb{Z}}$ finite (i.e. such that $\sum_{i\in \mathbb Z}\xi_i<\infty$).
\item[b)] Let $\mu$ be a probability measure on $\mathbb{R}_+^{\mathbb{Z}}$, then we say that it is a {\em  ``good infinite-volume  measure''} for ABEP$(\si,k)$ iff it concentrates on good infinite-volume configurations.
\end{itemize}
\ed
\bp\label{InfVol1}
\vskip.1cm
\noindent
\begin{itemize}
\item[1)] If $x \in \mathbb{R}_+^{\mathbb{Z}}$ is a { ``good infinite-volume configuration''} for ABEP$(\si,k)$ and  $\xi^{(\ell_1, \ldots, \ell_n)}$ is the configurations with
$n$ particles located at sites $\ell_1, \ldots, \ell_n\in \mathbb Z$, then the limit
\be
 \lim_{L\to\infty} \prod_{m=1}^n e^{2\sigma E_{\ell_m +1}(x^{(L)})} \; \mathbb{E}_{x^{(L)}}\left[D^\si(x(t),\xi^{(\ell_1, \ldots, \ell_n)})\right]
\ee
is well-defined for all $t\ge 0$ and is equal to
\begin{equation}
\lim_{L \to \infty} \prod_{m=1}^n e^{2\si E_{\ell_m +1}(x^{(L)})} \;\mathbb E_{\xi^{(\ell_1, \ldots, \ell_n)}} \left[D^\si(x^{(L)},\xi(t))\right]
\end{equation}
\item[2)] If $x \in \mathbb{R}_+^{\mathbb{Z}}$ is bounded, i.e.  $\sup_{i\in\mathbb{Z}} x_i <\infty$, then it is a ``good infinite-volume configuration'' for ABEP$(\si,k)$.
\item[3)] Let us denote by $\mathcal N_\lambda(t)$ a Poisson process of rate $\lambda>0$, and by $\mathbf E[\cdot]$ the expectation w.r. to its probability law. If $\mu$ is a probability measure on $\mathbb{R}_+^{\mathbb{Z}}$ such that for any $\lambda>0$ the expectation
\begin{equation}\label{mu1}
\mathbb{E}_\mu \left[ \mathbf E \left[ e^{\sum_{i=1}^{\mathcal N_\lambda(t)}x_{\ell +i}}\right]\right]
\end{equation}
is finite for all $t \ge 0$ and for any $\ell \in \mathbb Z$, then $\mu$ is a ``good infinite-volume  measure'' for ABEP$(\si,k)$.
\end{itemize}
\ep
\bpr
The proof is analogous to the proof of Proposition \ref{InfVol}.
\epr
\noindent
Later on  for a ``good'' infinite-volume configuration $x \in \mathbb{R}_+^{\mathbb{Z}}$ we will write
\begin{equation}\label{Notation}
\prod_{i \in \mathbb Z}e^{2\si\xi_i E_{i+1}(x)}\;\mathbb E_\xi \left[D^\si(x,\xi(t))\right]:=\lim_{L \to \infty} \prod_{i \in \mathbb Z}e^{2\si\xi_i E_{i+1}(x^{(L)})}\;\mathbb E_\xi \left[D^\si(x^{(L)},\xi(t))\right]
\end{equation}
and
\be\label{Notation1}
\prod_{m=1}^n e^{2\sigma E_{\ell_m +1}(x)} \; \mathbb{E}_{x}\left[D^\si(x(t),\xi^{(\ell_1, \ldots, \ell_n)})\right]:= \lim_{L\to\infty} \prod_{m=1}^n e^{2\sigma E_{\ell_m +1}(x^{(L)})} \; \mathbb{E}_{x^{(L)}}\left[D^\si(x(t),\xi^{(\ell_1, \ldots, \ell_n)})\right]
\ee

\subsection{$e^{-\si}$-exponential moment of the current of ABEP$(\si,k)$}
We start by defining the current for the ABEP$(\sigma,k)$ process on $\mathbb{Z}$.

\bd[Current]
 Let $\{x(t), \; t \ge 0\}$ be a c\`adl\`ag trajectory on the infinite-volume  configuration space $\mathbb R_+^{\mathbb Z}$, then the total integrated current $J_i(t)$ in the time interval $[0,t]$ is defined as total energy crossing the bond $(i-1,i)$
in the right direction.
\be \label{J2}
J_i(t)= E_i(x(t))-E_i(x(0)):= \lim_{L \to \infty}\(E_i(x^{(L)}(t))-E_i(x^{(L)})\)
\ee
where $E_i(x)$ is defined in \eqref{E} and  $x^{(L)}$ as in \eqref{xL}.
\ed

\bl[Current] \label{J00}
We have
$
\lim_{i\to -\infty} J_i(t) = 0$.
\el
\bpr
It immediately follows from the conservation of the total energy.
\epr

\bp[Current exponential moment via a dual walker]
\label{Lemma:E2}
The first exponential moment of $J_i(t)$ when the process is started from a ``good infinite-volume initial configuration''
$x \in \mathbb{R}_+^{\mathbb Z}$ at time $t=0$ is given by
\be
\mathbb{E}_{x}\[e^{-2\sigma J_i(x(t))}\] =  e^{-4kt} \,\sum_{n \in \mathbb{Z}} e^{-2\sigma (E_{n}(x)-E_i(x))} \;  I_{|n-i|}(4kt)
\label{Quiii2}
\ee
where  $I_n(t)$ is the modified Bessel function.

\ep

\bpr
Let  $\xi^{(\ell)} \in \mathbb R_+^{\mathbb Z}$  be  the configuration with a single  particle at site $\ell$.
Since the ABEP$(\sigma,k)$ is dual to the SIP$(2k)$ the dynamics of the single dual
particle is given by a continuous time symmetric random walker
$\ell(t)$ on $\mathbb Z$
jumping  at rate $2k$.
Since $x$ is a good configuration we have that the normalized expectation
$$
e^{2\sigma E_i(x)}\;\mathbb{E}_{x}\left[D(x(t),\xi^{(\ell)})\right] = \frac{1}{4k\sigma} \;e^{2\sigma E_i(x)}\;
\mathbb{E}_{x}\left[e^{-2\sigma E_{\ell+1}(x(t))}-e^{-2 \sigma E_{\ell}(x(t))}\right]
$$
and, from the duality relation  \eqref{one-dual} this is also equal to:
$$
e^{2\sigma E_i(x)}\;\mathbb{E}_{\xi^{(\ell)}}\left[D(x,\xi^{(\ell(t))})\right]= \frac{1}{4k\sigma} \; e^{2\sigma E_i(x)}\;\mathbf{E}_{\ell}\left[e^{-2 \sigma E_{\ell(t)+1}(x)}-e^{-2\sigma E_{\ell(t)}(x)}\right]
$$
where $\mathbf E_\ell$ denotes the expectation
with respect to the law of $\ell(t)$ started at site $\ell\in\mathbb{Z}$ at time $t=0$.
As a consequence,  for any $\ell \in \mathbb Z$
\begin{equation}\label{R2}
e^{2\sigma E_i(x)}\;\mathbb{E}_{x}\left[e^{-2 \sigma E_{\ell+1}(x(t))}\]=e^{2\sigma E_i(x)}\;\mathbb{E}_{x}\left[e^{-2 \sigma E_{\ell}(x(t))}\right]+e^{2\sigma E_i(x)}\; \mathbf{E}_{\ell}\left[e^{-2\sigma E_{\ell(t)+1}(x)}-e^{-2\sigma E_{\ell(t)}(x)}\right]
\end{equation}
from which it follows
\begin{eqnarray}
e^{2\sigma E_i(x)}\;\mathbb{E}_{x}\[e^{-2\sigma E_i(x(t))}\] &=&e^{2\sigma E_i(x)}\;\sum_{\ell\le i-1}\mathbf{E}_{\ell}\left[e^{-2\sigma E_{\ell(t)+1}(x)}-e^{-2\sigma E_{\ell(t)}(x)}\right] \nn \\
&=& e^{2\sigma E_i(x)}\;\sum_{\ell\le i-1}\mathbf{E}_{0}\left[e^{-2\sigma E_{\ell(t)+\ell+1}(x)}-e^{-2\sigma E_{\ell(t)+\ell}(x)}\right] \nn \\
&=& e^{2\sigma E_i(x)}\;\sum_{m\le i}\mathbf{E}_{0}\left[e^{-2\sigma E_{\ell(t)+m}(x)}\right]-\sum_{\ell\le i-1}\mathbf{E}_{0}\left[e^{-2\sigma E_{\ell(t)+\ell}(x)}\right] \nn \\
&=& e^{2\sigma E_i(x)}\; \mathbf{E}_{0}\left[e^{-2\sigma E_{\ell(t)+i}(x)}\right] \nn \\
&=& e^{2\sigma E_i(x)}\;\mathbf{E}_{i}\left[e^{-2\sigma E_{\ell(t)}(x)}\right]\;.
\label{Quii12}
\end{eqnarray}
Thus we have arrived to 
\be
\mathbb{E}_{x}\[e^{-2\sigma J_i(t)}\] =\mathbf{E}_{i}\left[e^{-2\sigma \(E_{\ell(t)}(x)-E_i(x)\)}\right]
\label{Quii2}
\ee
and the result \eqref{Quiii2} follows since
$$
\mathbf{E}_i(f(\ell(t)) = \sum_{n\in\mathbb{Z}} f(n) \cdot \mathbf{P}_i(\ell(t)=n)
$$
with
\begin{eqnarray}
\label{Bessel2}
\mathbf{P}_i(\ell(t) = n)
& = &
\mathbb{P}(\ell(t)=n \;|\; \ell(0) = i) \nonumber \\
& = & e^{-4kt}  I_{|n-i|}(4kt)
\end{eqnarray}
where  $I_n(x)$ is the modified Bessel function.
\epr

\br \label{Rem:LDP2}
Let $\ell(t)$ be a continuous time symmetric random walk
on $\mathbb Z$
jumping  at rate $2k$,  then \eqref{LDP} holds with
\be \label{I2}
{\cal I}(x)= 4k- \sqrt{x^2+(4k)^2} + x \log \left\{\frac{1}{4k}\[x+ \sqrt{x^2+(4k)^2}\]\right\}
\ee
\er

\noindent
We denote by  ${\mathbb E}^{\otimes \mu}$ the expectation of the ABEP$(\sigma,k)$ process on $\mathbb{Z}$
initialized with the  omogeneous product measure on $\mathbb{R}^{\mathbb Z}$ with marginals $\mu$
at time 0, i.e.
\be
{\mathbb E}^{\otimes \mu}[f(x(t))]= \int \( \otimes_{i\in\mathbb{Z}} \mu(dx_i) \) \; \mathbb E_x[f(x(t))]
\ee

\bp[Exponential moment for product initial condition]
\label{current-prod2}
Consider a probability measure $\mu$ on $\mathbb{R}^+$. Then,
 for the infinite volume ABEP$(\sigma,k)$, we have
\be \label{Exp-mom-prod2}
\mathbb{E}^{\otimes\mu}\[e^{-2\sigma J_i(t)}\]= \mathbf P_0 \[\ell(t)=0\]+
\mathbf E_0 \[\(\lambda_{+}^{\ell(t)}+ \lambda_-^{\ell(t)}\)\mathbf 1_{\ell(t)\ge 1}\]
\ee
where $\lambda_\pm:= \int  \mu(dy) e^{\pm 2 \sigma y}$
and $\ell(t)$ is the random walk defined in Remark \ref{Rem:LDP2}.
In particular we have
\be\label{M2}
\lim_{t \to\infty} \frac{1}{t} \log\mathbb E^{\otimes\mu}[e^{-2 \sigma J_i(t)}]= \sup_{x \ge 0}\left\{x \log \lambda_+ -{\cal I}(x)\right\} -\inf_{x \ge 0}{\cal I}(x)
\ee
with  ${\cal I}(x)$ given by \eqref{I2}.
\ep
\bpr It is easy to check that an homogeneous product measure $\mu$ verifies the condition \eqref{mu1} in Proposition \ref{InfVol}, thus it is a good infinite-volume probability measure for ABEP$(\si,k)$ in the sense of Definition \ref{Good1}. Thus  we can apply Proposition \ref{Lemma:E2}, in particular
from \eqref{Quii2} we have
\begin{eqnarray}
&&\mathbb{E}^{\otimes\mu}\[e^{-2\sigma J_i(t)}\]= \int \otimes\mu(dx)\, \mathbb E_x\[e^{-2\sigma J_i(t)}\]
\nn \\
&& = \int \otimes\mu(dx) \mathbf{E}_{i}\left[e^{-2\sigma \(E_{\ell(t)}(x)-E_i(x)\)}\right]=\nn \\
&&=\sum_{n \in \mathbb{Z}} \mathbf{P}_{i}\(\ell(t)=n\) \;\int \otimes\mu(dx)  e^{-2\sigma (E_{n}(x)-E_i(x))}\;. \nn \\
\end{eqnarray}
Since
\be
\int \otimes\mu(d \eta) e^{-2\sigma(E_x(\eta)-E_i(\eta))} = \lambda_-^{i-n} \; \mathbf 1_{\{n \le i\}} + \lambda_{+}^{n-i} \; \mathbf 1_{\{n>i\}}
\ee
it follows that
\begin{eqnarray}\label{Ct2}
\mathbb{E}^{\otimes\mu}\[e^{-2\sigma J_i(t)}\]&=&\sum_{n \le i} \mathbf P_i\(\ell(t)=n\) \lambda_-^{i-n}+ \sum_{n \ge i+1} \mathbf P_i\(\ell(t)=n\)  \lambda_+^{n-i} \nn \\
&=&  \mathbf{E}_{i}\left[\lambda_-^{i-\ell(t)} \mathbf 1_{\ell(t)\le i}+ \lambda_+^{\ell(t)-i} \mathbf 1_{\ell(t)\ge i+1}\right] \nn \\
&=&  \mathbf{E}_{0}\left[\lambda_-^{-\ell(t)} \mathbf 1_{\ell(t)\le 0}+ \lambda_+^{\ell(t)} \mathbf 1_{\ell(t)\ge 1}\right] \nn
\\
&=&  \mathbf{E}_{0}\left[\lambda_-^{\ell(t)} \mathbf 1_{\ell(t)\ge 0}+ \lambda_+^{\ell(t)} \mathbf 1_{\ell(t)\ge 1}\right]
\end{eqnarray}
where the last identity  follows from the symmetry of $\ell(t)$. Then \eqref{Exp-mom-prod2} is proved.

\vskip.2cm
\noindent
In order to prove \eqref{M2} we  rewrite \eqref{Exp-mom-prod2}
as
\begin{eqnarray} \label{Exp-mom-prod12}
\mathbb{E}^{\otimes\mu}\[e^{-2\sigma J_i(t)}\]
& = & \mathbf E_0 \[\lambda_+^{\ell(t)} \mathbf 1_{\ell(t)\ge 0}\]  \left(1+\mathcal E_1(t)+\mathcal E_2(t)\right)
\end{eqnarray}
with
\be
\mathcal E_1(t) :=\frac{\mathbf E_0 \[\(\lambda_+^{\ell(t)} +\lambda_-^{\ell(t)} \)\, \mathbf 1_{\ell(t)\ge 1}\]}
{\mathbf E_0 \[\lambda_+^{\ell(t)} \mathbf 1_{\ell(t)\ge 0}\]}, \qquad \mathcal E_2(t):=\frac{\mathbf P_0\(x(t)=0\)}{\mathbf E_0 \[\lambda_+^{\ell(t)} \mathbf 1_{\ell(t)\ge 0}\]}
\nn
\ee
where for  $i=1,2$ there exists $c_i>0$ such that
\be\label{E2}
\sup_{t\ge 0}|\mathcal  E_i(t)| \le c_i
\ee
This and the result of  Remark \ref{Rem:LDP2} conclude the proof of\eqref{M2}.
\epr

\section{Algebraic construction of ASIP$(q,k)$ and proof of the self-duality}
In this section we give the full proof of Theorem \ref{mainself}.
It follows closely  the lines of \cite{CGRS}  however the algebra and co-product 
are different.

\subsection{Algebraic structure and symmetries}

\subsubsection*{The quantum Lie algebra {${\mathcal{U}}_q(\mathfrak{su}(1,1))$}}
\label{quantum-algebra}

For $q\in(0,1)$ we consider the algebra with generators $K^{+}, K^{-}, K^{0}$ satisfying the commutation
relations

\begin{eqnarray}
\label{comm-suq2}
[K^+,K^-]=-[2K^0]_q, \qquad  [K^0,K^\pm]=\pm K^\pm\;,
\end{eqnarray}
where $[\cdot,\cdot]$ denotes the commutator, i.e. $[A,B] = AB-BA$, and
\begin{equation}
[2K^0]_q :=\frac{q^{2K^0}-q^{-2K^0}}{q-q^{-1}}\;.
\end{equation}
This is the quantum Lie algebra \red{${\mathcal{U}}_q(\mathfrak{su}(1,1))$}, that in the limit $q\to 1$ reduces
to the Lie algebra $\mathfrak{su}(1,1)$.  The
Casimir element is
\be
\label{casimir}
C = [K^0]_q[K^0-1]_q-  K^+ K^-
\ee
A standard representation   of the quantum Lie algebra ${\mathcal{U}}_q(\mathfrak{su}(1,1))$
is given by
\begin{equation}
\label{stand-repr}
\left\{
\begin{array}{lll}
{K}^+ |n\rangle &=& \sqrt{[\eta+2k]_q [\eta+1]_q}\;| n +1 \rangle
\\
{K}^- |n\rangle &=& \sqrt{[\eta]_q [\eta+2k-1]_q} \;| n-1 \rangle
\\
{K}^0 |n\rangle &=& (\eta+k) \;| n \rangle \;.
\end{array}
\right.
\end{equation}
$k \in \mathbb N$.
Here the collection of column vectors $|n\rangle$, with $n \in\mathbb N$, denote
the standard  orthonormal basis  with respect to the Euclidean scalar product,
i.e.  $|n\rangle = (0,\ldots,0,1,0,\ldots, 0)^T$ with the element $1$ in the $n^{\text{th}}$
position and with the symbol $^{T}$ denoting transposition.
Here and in the following, with abuse of notation, we use the same symbol for a linear
operator and the matrix associated to it in a given basis.
In the representation \eqref{stand-repr}
the ladder operators ${K}^+$ and ${K}^-$
are
one the adjoint of the other, namely
\be\label{transp}
({K}^+)^* = {K}^-
\ee
 and
the Casimir element is given by the diagonal matrix
$$
{C} |n\rangle =[ k]_q[k-1]_q |n\rangle\;.
$$

\noindent
We also observe that the ${\mathcal{U}}_q(\mathfrak{su}(1,1))$ commutation relations in \eqref{comm-suq2}
can be rewritten as follows
\begin{eqnarray}
\label{comm-new}
&& q^{K_0} K^+ = q \;  K^+ q^{K_0} \\
&& q^{K_0} K^-= q^{-1}\, K^-q^{K_0} \nonumber\\
&& [K^+,K^-]=-[2K^0]_q \nonumber
\end{eqnarray}

\subsubsection*{Co-product structure}
\label{cooooo}

A co-product for the quantum Lie algebra ${\mathcal{U}}_q(\mathfrak{su}(1,1))$ is defined as the map
$\Delta: {\mathcal{U}}_q(\mathfrak{su}(1,1))\to {\mathcal{U}}_q(\mathfrak{su}(1,1)) \otimes {\mathcal{U}}_q(\mathfrak{su}(1,1))$
\begin{eqnarray}
\label{co-product2}
\Delta(K^{\pm}) & = & K^{\pm} \otimes  q^{-K^0} + q^{K^0} \otimes K^{\pm}\;, \nonumber \\
\Delta(K^0) & = & K^0 \otimes 1 +  1\otimes K^0\;.
\end{eqnarray}
The co-product is an isomorphism for the quantum Lie algebra \red{$U_q(\mathfrak{sl}_2)$}, i.e.
\be
\label{coproduct}
[\Delta(K^+),\Delta(K^-)]=-[2\Delta(K^0)]_q, \qquad  [\Delta(K^0),\Delta(K^\pm)]=\pm \Delta(K^\pm)\;.
\ee
Moreover it can be easily checked that the co-product satisfies the co-associativity property
\be
\label{co-ass}
(\Delta\otimes 1) \Delta = (1\otimes \Delta) \Delta \;.
\ee
\noindent
Since we are interested in extended systems we will work with the
tensor product over copies of the ${\mathcal{U}}_q(\mathfrak{su}(1,1))$ quantum algebra.
We denote by $K_i^{+}, K_i^{-}, K_i^{0}$, with $i\in\mathbb{Z}$, the generators
of the $i^{th}$ copy. Obviously algebra elements of different copies commute.
As a consequence of \eqref{co-ass}, one can define iteratively
\red{$\Delta^{n}: {\mathcal{U}}_q(\mathfrak{su}(1,1)) \to {\mathcal{U}}_q(\mathfrak{su}(1,1))^{\otimes (n+1)}$},
i.e. higher power of $\Delta$, as follows: for $n=1$,  from  \eqref{co-product2} we have
\begin{eqnarray}
\label{cooo-product}
\Delta(K_i^{\pm}) & = & K_i^{\pm} \otimes  q^{-K_{i+1}^0} + q^{K_i^0} \otimes K_{i+1}^{\pm} \nonumber \\
\Delta(K_i^0) & = & K_i^0 \otimes 1 +  1 \otimes K_{i+1}^0\;,
\end{eqnarray}
for $n\ge 2$,
\begin{eqnarray}
\label{co-product-L}
\Delta^{n}(K_i^{\pm}) & = & \Delta^{n-1}(K_i^{\pm}) \otimes q^{-K^0_{n+i}}  +  q^{\Delta^{n-1}(K_i^0)} \otimes K_{n+i}^{\pm}  \nonumber\\
\Delta^{n}(K_i^0) & = & \Delta^{n-1}(K_i^0)  \otimes 1 +  \underbrace{1\otimes\ldots\otimes 1}_{n \text{ times}} \otimes {K_{n+i}^0}\;.
\end{eqnarray}

\subsubsection*{The quantum Hamiltonian}
\label{q-h}

Starting from the quantum Lie algebra ${\mathcal{U}}_q(\mathfrak{su}(1,1))$ 
and the co-product structure we would like to
construct a linear operator (called ``the quantum Hamiltonian'' in the following and denoted
by $H^{\phantom x}_{(L)}$ for a system of length $L$) with the following properties:
\begin{enumerate}
\item it is ${\mathcal{U}}_q(\mathfrak{su}(1,1))$ symmetric, i.e. it admits non-trivial symmetries constructed from
the generators of the quantum algebra; the non-trivial symmetries can then be used
to construct self-duality functions;
\item it can be associated to a continuos time Markov jump process,
i.e. there exists a representation given
by a matrix with non-negative out-of-diagonal elements (which can therefore
be interpreted as the rates of an interacting particle systems) and with zero sum
on each column.
\end{enumerate}

\vspace{0.2cm}
\noindent
A natural candidate for the quantum Hamiltonian operator is obtained by applying
the co-product to the Casimir operator $C$ in \eqref{casimir}.
Using the co-product definition \eqref{co-product2},  simple
algebraic manipulations   yield the following definition.
\bd[Quantum Hamiltonian]
\label{def-qh}
For every $L\in\mathbb{N}$,  $L\ge 2$, we consider the operator $H^{\phantom x}_{(L)}$ defined by
\be
\label{hami}
H^{\phantom x}_{(L)}
:=  \sum_{i=1}^{L-1} H^{i,i+1}_{(L)}
= \sum_{i=1}^{L-1} \left( h^{i,i+1}_{(L)} + c_{(L)} \right) \;,
\ee
where the two-site Hamiltonian is the sum of
\be
\label{const}
c_{(L)} = \frac{(q^{2k}-q^{-2k})(q^{2k-1}-q^{-(2k-1)})}{(q-q^{-1})^2} \underbrace{1\otimes\cdots \otimes 1}_{L \text{ times}}
\ee
and
\be
h^{i,i+1}_{(L)} := \underbrace{1\otimes\cdots \otimes 1}_{(i-1) \text{ times}} \otimes \Delta(C_i) \otimes
\underbrace{1 \otimes \cdots \otimes 1}_{(L-i-1) \text{ times}}
\ee
and, from \eqref{casimir} and \eqref{co-product2},
\be
\Delta(C_i) = \Delta(K_i^+)\Delta(K_i^-)- \Delta([K_i^0]_q) \Delta([K_{i}^0-1]_q)
\ee
Explicitely
\begin{eqnarray}
\label{deltaci}
\Delta(C_i)
& =  &
 q^{K_i^0}\Bigg \{ K_i^+ \otimes K_{i+1}^-+ K_i^- \otimes K_{i+1}^+\Bigg \} q^{-K_{i+1}^0}+
 K_i^+ K_i^- \otimes q^{-2K^0_{i+1}} + q^{2K^0_{i}} \otimes K_{i+1}^+ K_{i+1}^-  \nonumber \\ & - & \frac{1}{(q-q^{-1})^2}
 \left\{q^{-1} q^{2K^0_i} \otimes q^{2K^0_{i+1}} +q   q^{-2K^0_i} \otimes q^{-2K^0_{i+1}} -(q+q^{-1}) \right\}
\quad\qquad
\end{eqnarray}
\ed
\br
By specializing \eqref{deltaci} to the representation \eqref{stand-repr} we get
\begin{eqnarray}
\label{deltaci11}
\Delta(C_i)
& =  &
q^{K_i^0}\Bigg \{ K_i^+ \otimes K_{i+1}^-+ K_i^- \otimes K_{i+1}^+ \: +\\
&& \hskip1.1cm -\: \frac {(q^k+q^{-k})(q^{k-1}+q^{-(k-1)})}{2(q-q^{-1})^{2}} \:
\(q^{K_i^0}-q^{-K_i^0}\) \otimes
\(q^{K_{i+1}^0}-q^{-K_{i+1}^0}\) \nonumber \\
& &
\quad\qquad - \: \frac{(q^k-q^{-k})(q^{k-1}-q^{-(k-1)})}{2(q-q^{-1})^{2}} \(q^{K_i^0}+q^{-K_i^0}\)\otimes \(q^{K_{i+1}^0}+q^{-K_{i+1}^0}\) \Bigg \} q^{-K_{i+1}^0} \nn
\end{eqnarray}
\er

\br The diagonal operator $c_{(L)}$ in \eqref{const} has been added so that the ground state $|0\rangle_{(L)} := \otimes_{i=1}^L |0\rangle_i$ is a
right eigenvector with eigenvalue zero,
i.e. $H_{(L)} |0\rangle_{(L)} = 0$ as it is immediately seen using \eqref{stand-repr}.
\er

\bp
\label{h-herm}
In the representation \eqref{stand-repr} the operator $H_{(L)}$ is self-adjoint.
\ep
\bpr
It is enough to consider the non-diagonal part of $H_{(L)}$. Using \eqref{transp} we have
\begin{eqnarray}
&&\(q^{K_i^0} K_i^+ \otimes K_{i+1}^- q^{-K_{i+1}^0} + q^{K_i^0} K_i^- \otimes K_{i+1}^+  q^{-K_{i+1}^0} \)^* \nn \\
&=& K_i^- q^{K_i^0} \otimes  q^{-K_{i+1}^0}  K_{i+1}^++  K_i^+ q^{K_i^0} \otimes q^{-K_{i+1}^0} K_{i+1}^- \nn \\
&=&  q^{K_i^0+1} K_i^- \otimes   K_{i+1}^+ q^{-K_{i+1}^0-1}  +   q^{K_i^0-1}  K_i^+ \otimes K_{i+1}^- q^{-K_{i+1}^0+1}  \nn
\end{eqnarray}
where the last identity follows by using the commutation relations  \eqref{comm-new}. This concludes the proof.
\epr

\subsubsection*{Basic symmetries}
\label{basic}

It is easy to construct symmetries for the operator $H^{\phantom x}_{(L)}$ by using the property that the co-product is
an isomorphism for the ${\mathcal{U}}_q(\mathfrak{su}(1,1))$ algebra.

\bt[Symmetries of $H_{(L)}$]
\label{theo-symm}
Recalling \eqref{co-product-L}, we define the operators
\begin{eqnarray}
K_{(L)}^{\pm} & := & \Delta^{L-1}(K_1^{\pm}) = \sum_{i=1}^L q^{K_1^0} \otimes \cdots \otimes q^{K_{i-1}^0} \otimes K_i^{\pm} \otimes q^{-K_{i+1}^0} \otimes \ldots \otimes q^{-K_L^0}\;,
\nonumber \\
K_{(L)}^{0} & := & \Delta^{L-1}(K_1^0) = \sum_{i=1}^L
 \underbrace{1\otimes\cdots \otimes 1}_{(i-1) \text{ times}} \otimes K_i^0 \otimes
\underbrace{1 \otimes \cdots \otimes 1}_{(L-i) \text{ times}}
\;.
\end{eqnarray}
They are symmetries of the Hamiltonian \eqref{hami}, i.e.
\be
[H_{(L)}^{\phantom x},K_{(L)}^{\pm}]= [H_{(L)}^{\phantom x},K_{(L)}^{0}] = 0\;.
\ee
\et
\begin{proof}
We proceed by induction and prove only the result for  $K_{(L)}^{\pm}$ (the case $K^{0}_{(L)}$ is similar).
By construction $K_{(2)}^{\pm} := \Delta(K^{\pm})$ are symmetries
of the two-site Hamiltonian $H^{\phantom x}_{(2)}$. Indeed this is an immediate consequence of the fact that
the co-product defined in \eqref{coproduct} conserves the commutation relations and the Casimir operator
\eqref{casimir} commutes with any other operator in the algebra :
$$
[H^{\phantom x}_{(2)}, K^{\pm}_{(2)}] = [\Delta(C_1),\Delta(K_1^{\pm})] = \Delta ( [C_1,K_1^{\pm}] ) = 0\;.
$$
For the induction step assume now that it holds $[H^{\phantom x}_{(L-1)}, K^{\pm}_{(L-1)}] = 0$. We have
\be
\label{eq-zero}
[H^{\phantom x}_{(L)}, K^\pm_{(L)}]  = [H^{\phantom x}_{(L-1)}, K^\pm_{(L)}]  + [h^{L-1,L}_{(L)}, K^\pm_{(L)}]
\ee
The first term on the right hand side of \eqref{eq-zero}  can be seen to be zero using \eqref{co-product-L} with $i=1$ and
$n=L-1$:
$$
[H^{\phantom x}_{(L-1)}, K^{\pm}_{(L)}] = [H^{\phantom x}_{(L-1)}, K^{\pm}_{(L-1)} q^{-K^0_{L}} + q^{K^0_{(L-1)}} K^\pm_{L}]
$$
Distributing the commutator with the rule $[A,BC] = B[A,C] + [A,B] C$, the induction hypothesis and the fact that spins on
different sites commute imply the claim.
The second term on the right hand side of \eqref{eq-zero} is also seen to be zero by writing
$$
[h^{L-1,L}_{(L)}, K^{\pm}_{(L)}] = [h^{L-1,L}_{(L)}, K^{\pm}_{(L-2)} q^{-\Delta(K^0_{L-1})} + q^{K^0_{(L-2)}} \Delta(K^\pm_{L-1})] = 0\;.
$$
\end{proof}

\br In the case $q=1$, the quantum Hamiltonian in Definition \ref{def-qh} reduces to
the (negative of the) well-known Heisenberg ferromagnetic quantum spin chain 
\begin{equation}
H_{(L)} = \sum_{i=1}^{L-1} \( K_i^+ K_{i+1}^-+ K_i^- K_{i+1}^+ -2 K^0_i K^0_{i+1}+2k^2 \)\;,
\end{equation}
with 
spins $K_i$ satisfying the \red{$\mathfrak{su}(1,1)$} Lie algebra.
The symmetries of this Hamiltonian are given by
$$
K^{\pm}_{(L)} = \sum_{i=1}^L K_i^{\pm} \qquad \text{and} \qquad K_{(L)}^{o} = \sum_{i=1}^L K_i^o\;.
$$
\er

\subsection{Construction of  ASIP$(q,k)$ from the quantum Hamiltonian}
\label{proc}

In order to construct a Markov process from the quantum Hamiltonian $H_{(L)}$, we make use of the following Theorem which has been proven in \cite{CGRS}.
\bt[Positive ground state transformation]\label{corooo}
Let $A$ be a {$|\Omega|\times |\Omega|$} matrix with non-negative off diagonal elements.
Suppose there exists a column vector
$e^{\psi}:= g \in \R^{{|\Omega|}}$
 with strictly positive entries and such that
$Ag=0$. Let us denote by $G$ the diagonal matrix with entries
{$G(x,x)=g(x)$ for $x\in\Omega$}.
Then we have the following
\begin{itemize}
\item[a)] The matrix
$$
\caL= G^{-1} A G
$$
with entries
{
\be\label{amodi}
\caL(x,y) = \frac{A(x,y) g(y)}{g(x)}, \qquad x,y \in \Omega\times\Omega
\ee
}
is the generator of a Markov process $\{X_t:t\geq 0\}$ taking values on {$\Omega$}.
\item[b)] $S$ commutes with $A$ if and only if $G^{-1} S G$ commutes with $\caL$.
\item[c)] If $A=A^*$, where $^*$ denotes transposition, then
the probability measure {$\mu$} on {$\Omega$}
{
\be\label{revmes}
\mu(x)= \frac{(g(x))^2}{\sum_{x\in\Omega} (g(x))^2}
\ee
}
is reversible for the process with generator $\caL$.

\end{itemize}
\et
\noindent

\vskip.5cm
\noindent
 Now we apply item a) of Theorem \ref{corooo} with $A=H_{(L)}$.
At this aim we  need a non-trivial symmetry which yields a non-trivial ground state.
Starting from the basic symmetries of $H_{(L)}$ described in Section
\ref{basic}, and inspired by the analysis of the symmetric case ($q=1$),
it will be convenient to consider the {\em exponential} of those symmetries.

\subsubsection*{The $q$-exponential and its pseudo-factorization}

\bd[$q$-exponential]
We define the $q$-analog of the exponential function  as
\begin{equation}
\label{q-exp}
{\exp}_q(x):= \sum_{n\ge 0} \frac{x^n}{\{n\}_q!}
\end{equation}
where
\begin{equation}
\label{q-num2}
\{n\}_q:= \frac{1-q^n}{1-q}
\end{equation}
\ed
\br The $q$-numbers in \eqref{q-num2} are related to the $q$-numbers in \eqref{q-num} by
the relation $\{n\}_{q^2}= [n]_q q^{n-1}$. This implies  $\{n\}_{q^2} ! = [n]_q! \, q^{n (n-1)/2}$
and therefore
\begin{equation}\label{exptilde}
{\exp}_{q^2}(x)=\sum_{n \ge 0} \frac{x^n}{[n]_q!}\, q^{-n(n-1)/2}
\end{equation}
\er
\vskip.1cm
\noindent
\bp[Pseudo-factorization]\label{lemma:E}
\label{pseudo}
Let $\{g_1,\ldots,g_L\}$ and $\{k_1,\ldots,k_L\}$ be  operators  such that for $L \in \mathbb N$
and $g\in\mathbb{R}$
\begin{equation}\label{r}
k_i g_i=r g_ik_i  \qquad\text{for}\quad i\in \Lambda_L\;.
\end{equation}
Define
\begin{equation}
g^{(L)}:=\sum_{i=1}^L  k^{(i-1)} g_i, \quad \text{with} \quad   k^{(i)}:= k_1 \cdot \dots \cdot k_{i} \quad   \text{for $i\ge 1$ and} \quad k^{(0)}=1,
\end{equation}
 then
\begin{equation}\label{fact1}
{\exp}_{r}(g^{(L)})=  {\exp}_{r}(g_1)\cdot {\exp}_{r}(k^{(1)} g_2) \cdot \dots \cdot {\exp}_{r}(k^{(L-1)} g_L)
\end{equation}
Moreover let
\begin{equation}
\hat g^{(L)}:=\sum_{i=1}^L g_i \, h^{(i+1)}, \quad \text{with} \quad   h^{(i)}:= k_i^{-1} \cdot \dots \cdot k^{-1}_{L} \quad   \text{for $i\le L$ and } \quad h^{(L+1)}=1,
\end{equation}
 then
\begin{equation}\label{fact2}
{\exp}_{r}(\hat g^{(L)})=  {\exp}_{r}(g_1 \, h^{(2)}) \cdot \dots \cdot {\exp}_{r}(g_{L-1} \, h^{(L)}) \cdot {\exp}_{r}(g_{L})
\end{equation}
\ep
\vskip.5cm
\noindent
See \cite{CGRS} for the proof.


\subsubsection*{The exponential symmetry $S^+_{(L)}$}
\label{exps+}

In this Section we identify the symmetry that will be used in the construction of the process ASIP$(q,k)$.
To have a symmetry that has quasi-product form over the sites we preliminary define more convenient
generators of the ${\mathcal{U}}_q(\mathfrak{su}(1,1))$ quantum Lie algebra.
Let
\begin{equation}\label{EF2}
 E := q^{K^0} K^+ , \qquad  F:= K^-  q^{-K^0}\qquad \text{and} \qquad  K:= q^{2K^0}
\end{equation}
From the commutation relations \eqref{comm-suq2} we deduce that $(E, F, K)$  verify the relations
\begin{equation}\label{commRel2}
KE= q^2E K \qquad \text{and} \qquad K F = q^{-2} F K  \qquad [E,F] = - \,\frac{K-K^{-1}}{q-q^{-1}}\;.
\end{equation}
Moreover, from Theorem \ref{theo-symm}, the following co-products
\begin{eqnarray}
&&\Delta(E_1):=\Delta(q^{K^0_1}) \cdot \Delta(K^+_1)  = E_1 \otimes  \mathbf 1 + K_1 \otimes   E_2\label{E2}\\
&&\Delta(F_1):= \Delta(K^-_1) \cdot \Delta(q^{-K^0_1})  = F_1 \otimes  K_2^{-1} + \mathbf 1 \otimes   F_2\label{F2}
\end{eqnarray}
are still symmetries of $H_{(2)}$.
In general we can extend  \eqref{E2} and \eqref{F2} to $L$ sites, then we have that
\begin{eqnarray}\label{EL2}
 E^{(L)}&:=& \Delta^{(L-1)} (E_1) \nn \\
&=& \Delta^{(L-1)}(q^{K^0_1}) \cdot \Delta^{(L-1)}(K^+_1) \nn \\
&=& q^{K^0_1} K^+_1 + q^{2K^0_1 + K^0_2} K^+_2 + ...+ q^{2 \sum_{i=1}^{L-1} K^0_i + K^0_L} K^+_L \nn \\
&=&  E_1 +K_1   E_2 + K_1 K_2  E_3 + ...+ K_1 \cdot ... \cdot K_{L-1}  E_L
\end{eqnarray}
\begin{eqnarray}\label{FL234}
 F^{(L)}&:=& \Delta^{(L-1)} (F_1) \nn \\
&=&  \Delta^{(L-1)}(K^-_1)  \cdot \Delta^{(L-1)}(q^{-K^0_1})\nn \\
&=& K^-_1 q^{-K^0_1-2\sum_{i=2}^L K_i^0} + \dots + K^-_{L-1} q^{-K^0_{L-1}-2K^0_L}  + K^-_L q^{-K^0_L}  \nn \\
&=&  F_1 \cdot K_2^{-1} \cdot ... \cdot K_{L}^{-1}+\dots +  F_{L-1}\cdot K_L^{-1}+   F_L
\end{eqnarray}
are symmetries of $H$. If we consider now the symmetry obtained by $q$-exponentiating $E^{(L)}$ then
this operator will pseudo-factorize by Proposition \ref{pseudo}.
\begin{lemma}
\label{lemma:S+}
 The operator
 \begin{equation}
 \label{s+}
 S_{(L)}^+:={\exp}_{q^{2}}( E^{(L)})
\end{equation}
is a  symmetry of $H_{(L)}$.
Its matrix elements are given by
\begin{equation}\label{Selem}
\langle \eta_1, ...,\eta_L |  S^+_{(L)} | \xi_1, ...,\xi_L \rangle=
\prod_{i=1}^L \sqrt{\binom{\eta_i}{\xi_i}_q \binom{\eta_i+2k-1}{\xi_i+2k-1}_q} \cdot \mathbf 1_{\eta_i \ge \xi_i}\, q^{(\eta_i-\xi_i)\left[1 +k + \xi_i +2\sum_{m=1}^{i-1}(\xi_m+k) \right]}
\end{equation}
\end{lemma}
\vskip.5cm
\noindent
\bpr
From \eqref{commRel2} we know that the operators   $ E_i, K_i$, copies of the operators defined in \eqref{EF2}, verify the conditions \eqref{r} with $r=q^{2}$.
As a consequence, from \eqref{EL2}, \eqref{s+} and Proposition \ref{lemma:E}, we have
\begin{eqnarray}
 S^+_{(L)}
 &=&{\exp}_{q^{2}}( E^{(L)})\nn \\
&=& {\exp}_{q^{2}}( E_1)\cdot {\exp}_{q^{2}}(K_1  E_2) \cdots {\exp}_{q^{2}}(K_1 \cdots K_{L-1}  E_L) \nn\\
&=&{\exp}_{q^{2}} \(q^{K^0_1} K^+_1\) \cdot {\exp}_{q^{2}}\(q^{2K^0_1} q^{K^0_2} K^+_2\) \cdots   {\exp}_{q^{2}}\(q^{2 \sum_{i=1}^{L-1} K^0_i + K^0_L} K^+_L\) \nn \\
&=&  S_1^+  S_2^+ \cdots  S_L^+
\end{eqnarray}
\vskip.1cm
\noindent
where $S^+_i:=  {\exp}_{q^{2}}\(q^{2 \sum_{m=1}^{i-1} K^0_m + K^0_i} K^+_i\)$ has been defined. Using \eqref{exptilde}, we find
\begin{eqnarray}
 S^+_i | \xi_1,\ldots,\xi_L \rangle &=&   \sum_{\ell_i \ge 0} \frac{1}{[\ell_i]_q!}\(q^{2 \sum_{m=1}^{i-1}K_m^0 + K_i^0}K_i^+\)^{\ell_i} q^{-\frac 1 2 \ell_i(\ell_i-1)} | \xi_1,\ldots ,\xi_L\rangle \\
&& \hskip-2.5cm= \sum_{\ell_i \ge 0} \sqrt{\binom{\xi_i+\ell_i}{\ell_i}_q\cdot \binom{\xi_i+2k+\ell_i-1}{\ell_i}_q} \cdot  q^{\ell_i(\xi_i+k+1)+ 2 \ell_i \sum_{m=1}^{i-1}(\xi_m+k)} |\xi_1,\ldots, \xi_i+\ell_i, \ldots, \xi_L \rangle \nn
\end{eqnarray}
where in the last equality we used \eqref{stand-repr}. Thus we find
\begin{eqnarray}
\label{action-essepiuuu}
 S^+_{(L)} | \xi_1,\ldots,\xi_L \rangle &=&   S^+_1  S_2^+ \dots  S_L^+ | \xi_1,\dots,\xi_L \rangle  \\
&=&  \sum_{\ell_1, \ell_2, \dots, \ell_L \ge 0}\prod_{i=1}^L \Bigg(\sqrt{\binom{\xi_i+\ell_i}{\ell_i}_q\cdot \binom{\xi_i+2k+\ell_i-1}{\ell_i}_q} \nn \\
&&  \hskip2cm\cdot  \; q^{\ell_i(\xi_i+k+1)+ 2 \ell_i \sum_{m=1}^{i-1}(\xi_m+k)}  \Bigg) |\xi_1+\ell_1,\ldots, \xi_L+\ell_L \rangle  \nn
\end{eqnarray}
form which the matrix elements in \eqref{Selem} are immediately found.
\epr

\subsubsection*{Construction of a positive ground state and  the associated Markov process ASEP$(q,j)$}
\label{sec5.3}

By applying Theorem \ref{corooo} we are now ready to identify the
stochastic process related to the Hamiltonian $H_{(L)}$ in \eqref{hami}.

\noindent
We start from the  state ${\bf |0\rangle} = |0,\ldots,0\rangle$ which
is obviously a trivial ground state of $H_{(L)}$.
We then produce a positive ground state by acting with
the symmetry ${S}^{+}_{(L)}$ in \eqref{s+}.
Using \eqref{action-essepiuuu} we obtain
\begin{eqnarray}\label{g}
|g\rangle &=&  S^+_{(L)} |0, \ldots,0\rangle = \sum_{\ell_1, \ell_2, \dots, \ell_L \ge 0}\prod_{i=1}^L \sqrt{\binom{2k+\ell_i-1}{\ell_i}_q}
\cdot  \; q^{\ell_i(1-k+ 2 k i)} \; |\ell_1, ..., \ell_L \rangle  \nonumber
\end{eqnarray}
Following the scheme in Theorem \ref{corooo} we construct the operator $G_{(L)}$ defined by
\be
G_{(L)} |\eta_1,\ldots,\eta_L \rangle =  |\eta_1,\ldots,\eta_L\rangle \langle \eta_1,\ldots,\eta_L | {S}^+ | 0,\ldots,0 \rangle
\ee
In other words $G_{(L)}$ is represented by a diagonal matrix whose coefficients in the standard basis read
\begin{equation}\label{G}
\langle \eta_1,\ldots,\eta_L| G_{(L)} | \xi_1,\ldots,\xi_L \rangle= \prod_{i=1}^L\sqrt{\binom{\eta_i+2k-1}{\eta_i}_q} \cdot  q^{\eta_i(1-k+2k i)} \cdot \delta_{\eta_i=\xi_i}
\end{equation}
Note that $G_{(L)}$ is factorized over the sites, i.e.
\be
\langle \eta_1,\ldots,\eta_L| G_{(L)} | \xi_1,\ldots,\xi_L \rangle= \otimes_{i=1}^L \langle \eta_i| G_i| \xi_i\rangle
\ee
As a consequence of item a) of Theorem \ref{corooo},
the operator ${\cal L}^{(L)}$ conjugated to $H_{(L)}$ via $G^{-1}_{(L)}$, i.e.
\be
\label{tildehami}
{\cal L}^{(L)} = G^{-1}_{(L)} H_{(L)} G_{(L)}
\ee
is the generator of a Markov jump process $\eta(t) = (\eta_1(t),\ldots,\eta_L(t))$
describing particles jumping on the
chain $\Lambda_L$. The state space of such a process is given by $\Omega_L$
and its elements are denoted by  $\eta =  (\eta_1,\ldots,\eta_L)$, where
$\eta_i$  is interpreted as the number of particles at site $i$.
The asymmetry is controlled by the parameter
$0 < q \le 1$.
\vskip.4cm

\bp
The action of the Markov generator  ${\cal L}^{(L)}:=  G^{-1}_{(L)} H_{(L)} G_{(L)}$ is given by
\begin{eqnarray}
\label{gennn}
&&({\cal L}^{(L)}f)(\eta) =\sum_{i=1}^{{L-1}} ({\cal L}_{i,i+1}f)(\eta) \qquad \text{with} \nonumber\\
({\cal L}_{i,i+1}f)(\eta)
& = & q^{\eta_i-\eta_{i+1}+(2k-1)} [\eta_i]_q [2k+\eta_{i+1}]_q (f(\eta^{i,i+1}) - f(\eta)) \nonumber \\
& +  & q^{\eta_i-\eta_{i+1}-(2k-1)} [2k+\eta_i]_q [\eta_{i+1}]_q (f(\eta^{i+1,i}) - f(\eta))
\end{eqnarray}
\ep
\bpr
From  Proposition \ref{h-herm} we know that $H_{(L)}^*=H_{(L)}$, hence we have that the operator $\tilde{H}_{(L)}:=G_{(L)}H_{(L)}G^{-1}_{(L)}$  is the transposed of the generator $\caL^{(L)}$ defined by \eqref{tildehami}. Then we have to verify that the transition rates to move from $\eta$ to $\xi$ for the Markov process generated by \eqref{gennn} are equal to the elements $\langle \xi|\tilde{H}_{(L)}|\eta \rangle$.

\noindent
Since we already know that $\caL^{(L)}$ is a Markov generator, in order to prove the result
it is sufficient to apply the similarity transformation given by the matrix $G_{(L)}$ defined in \eqref{G}
to the non-diagonal terms of \eqref{deltaci}, i.e. $q^{K^0_{i}} K_i^\pm K_{i+1}^\mp q^{-K^0_{i+1}}$.  We show here the computation only for the first term, being the computation for the other term  similar. \\
 We have
\begin{eqnarray}
& &
\langle \xi_i,\xi_{i+1} | G_i G_{i+1} \cdot q^{K^0_{i}} K_i^+ K_{i+1}^-q^{-K^0_{i+1}} \cdot G^{-1}_i G^{-1}_{i+1}  | \eta_i,\eta_{i+1}\rangle \nn\\
& & =
\langle \xi_i | G_i q^{K^0_{i}} K_i^+ G^{-1}_i  | \eta_i\rangle \otimes
\langle \xi_{i+1} | G_{i+1} K_{i+1}^-q^{-K^0_{i+1}} G^{-1}_{i+1}  | \eta_{i+1}\rangle
\end{eqnarray}
Using \eqref{G} and \eqref{stand-repr} one has
\begin{eqnarray}
\langle \xi_i | G_i q^{K^0_{i}} K_i^+ G^{-1}_i  | \eta_i\rangle  =q^{\eta_i +2 +2ki}\; [2k+\eta_i]_q  \langle \xi_i  | \eta_i +1 \rangle
\end{eqnarray}
and
\begin{eqnarray}
\langle \xi_{i+1} | G_{i+1} K_{i+1}^-q^{-K^0_{i+1}} G^{-1}_{i+1}  | \eta_{i+1}\rangle =
q^{-\eta_{i+1} -2k-1 -2ki}\; [\eta_{i+1}]_q  \langle \xi_{i+1}  | \eta_{i+1} -1 \rangle
\end{eqnarray}
Multiplying the last two expressions one has
\be
\langle \eta^{i+1,i}|  {{\tilde{H}_{(L)}}}|\eta\rangle = q^{\eta_i-\eta_{i+1}-2k+1} [2k+\eta_i]_q [\eta_{i+1}]_q
\ee
that corresponds indeed to the rate to move  from $\eta$ to $\eta^{i+1,i}$ in \eqref{gennn}. This concludes the proof.
\epr
\vskip.3cm

\br \label{natural}
From item c) of Theorem \ref{corooo}, we have that  the product measure $\mu_{(L)}$ defined by
\be
\mu_{(L)}(\eta)=  \langle \eta|G_{(L)}^2|\eta \rangle
\ee
is a reversible  measure of ${\cal L}^{(L)}$. Notice that it corresponds to the reversible measure $\mathbb P^{(\alpha)}$ defined in \eqref{stat-meas} with the choice $\alpha=1$.

\er

\subsection{Self-Duality of  ASIP$(q,k)$}
\label{D}

%
%
The following Proposition  has been proven in \cite{CGRS} and it will be key to the
proof of ASIP$(q,k$ self-duality.
\bp\label{dualprop}
Let $A=A^*$ be a matrix with non-negative off-diagonal elements, and $g$
an eigenvector of $A$ with eigenvalue zero, with strictly positive entries. Let $\caL= {G^{-1}AG}$ be the corresponding
Markov generator.
Let $S$ be a symmetry of $A$, then $G^{-1}SG^{-1}$  is a self-duality function for the process with generator $\caL$.
\ep

We now use Proposition \ref{dualprop} and the exponential simmetry  obtained in Section \ref{exps+} to deduce a non-trivial duality function for the ASIP$(q,k)$ process.

%
%
%
%
%
%
\vskip.3cm
\noindent
{\bf{\small{P}{\scriptsize ROOF OF \eqref{dualll} IN THEOREM} \ref{mainself}}}.
From Proposition \ref{h-herm} we know that $H_{(L)}$ is self-adjoint, then, using Proposition \ref{dualprop} with $A=H_{(L)}$, $G=G_{(L)}$ given by \eqref{G} and $S=S_{(L)}^+$  given by \eqref{Selem} it follows that
\be
G^{-1}_{(L)} S^{+}_{(L)} G^{-1}_{(L)}
\ee
is a self-duality function for the process generated by ${\cal L}^{(L)}$.
 Its elements are computed as follows:
\begin{eqnarray}
&& \langle \eta |G_{(L)}^{-1} S_{(L)}^+ G_{(L)}^{-1}| \xi \rangle =\\
&& = \prod_{i=1}^L\(\sqrt{\binom{2k+\eta_i-1}{\eta_i}_q} \cdot  q^{\eta_i(1-k+ 2 k i)} \)^{-1} \langle \eta |  S^+_i | \xi \rangle  \(\sqrt{\binom{2k+\xi_i-1}{\xi_i}_q} \cdot  q^{\xi_i(1-k+ 2 k i)}\)^{-1}= \nn\\
&& = \prod_{i=1}^L   \sqrt{\binom{\eta_i}{\xi_i}_q \binom{\eta_i+2k-1}{\xi_i+2k-1}_q \bigg/  \binom{2k+\eta_i-1}{\eta_i}_q \binom{2k+\xi_i-1}{\xi_i}_q }  \cdot \nn \\
&& \hskip5cm \cdot \; q^{(\eta_i-\xi_i)\left[2\sum_{m=1}^{i-1}(\xi_m+k) +\xi_i+k+1\right]- (2 k i -k +1) (\eta_i+\xi_i)} \cdot \mathbf 1_{\xi_i \le \eta_i} =\nonumber \\
&& = q^{2\sum_{i=1}^L (k\xi_i-\eta_i)}\;\prod_{i=1}^L   {\frac{[2k-1]_q![\eta_i]_q!}{[\xi_i+2k-1]_q![\eta_i-\xi_i]_q!}}  \cdot \; q^{(\eta_i-\xi_i)\left[2\sum_{m=1}^{i-1}\xi_m +\xi_i\right]- 4 k i \xi_i} \cdot \mathbf 1_{\xi_i \le \eta_i} \nonumber
\end{eqnarray}
Since both the original process and the dual process conserve the total number of particles it follows that $D_{(L)}$ in \eqref{dualll}
is also a duality function. \qed

\end{document}